\documentclass[12pt]{article}
\begin{document}

\newtheorem{example}{Example}[section]
\newtheorem{theorem}{Theorem}[section]
\newtheorem{lemma}{Lemma}[section]
\newtheorem{corollary}{Corollary}[section]
\newtheorem{proposition}{Proposition}[section]
\newtheorem{remark}{Remark}[section]
\setlength{\textwidth}{13cm}
\def \halmos{\hfill\mbox{ qed}\\}
\newcommand{\eqnsection}
{\renewcommand{\theequation}{\thesection.\arabic{equation}}
\makeatletter \csname  @addtoreset\endcsname{equation}{section}
\makeatother}

\def \njc{{\bf !!! note Jay's change  !!!}}
\def  \enc{{\bf end of current changes }}
\def \nnc{{\bf !!! note NEW	 change  !!!}}
\def\square{{\vcenter{\vbox{\hrule height.3pt
                      \hbox{\vrule width.3pt height5pt \kern5pt
                         \vrule width.3pt}
                      \hrule height.3pt}}}}
\def \grad{\bigtriangledown}
\def \nc{{\bf !!! note change !!! }}
          \def \Proof{\noindent{\bf Proof}$\quad$}
\newcommand{\re}[1]{(\ref{#1})}
\def \ov{\overline}
\def \un{\underline}
\def \be{\begin{equation}}
\def \ee{\end{equation}}
\def \bt{\begin{theorem}}
\def \et{\end{theorem}}
\def \bc{\begin{corollary}}
\def \ec{\end{corollary}}
\def \br{\begin{remark} }
\def \er{ \end{remark}}
\def \bl{\begin{lemma}}
\def \el{\end{lemma}}
\def \bex{\begin{example}}
\def \eex{\end{example}}
\def \bea{\begin{eqnarray}}
\def \eea{\end{eqnarray}}
\def \bas{\begin{eqnarray*}}
\def \eas{\end{eqnarray*}}
\def \al{\alpha}
\def \bb{\beta}
\def \ga{\gamma}
\def \Ga{\Gamma}
\def \de{\delta}
\def \De{\Delta}
\def \ep{\epsilon}
\def \vep{\varepsilon}
\def \la{\lambda}
\def \La{\Lambda}
\def \ka{\kappa}
\def \om{\omega}
\def \Om{\Omega}
\def \va{\varrho}
\def \ffi{\Phi}
\def \vf{\varphi}
\def \si{\sigma}
\def \Si{\Sigma}
\def \vsi{\varsigma}
\def \th{\theta}
\def \Th{\Theta}
\def \ups{\Upsilon}
\def \ze{\zeta}
\def \tr{\nabla}
\def \ff{\infty}
\def \wh{\widehat}
\def \wt{\widetilde}
\def \dar{\downarrow}
\def \rar{\rightarrow}
\def \uar{\uparrow}
\def \sbs{\subseteq}
\def \mpt{\mapsto}
\def \R{{\bf R}}
\def \G{{\bf G}}
\def \H{{\bf H}}
\def \Z{{\bf Z}}
\def \S{{\bf S}}
\def \sfB{{\sf B}}
\def \sfS{{\sf S}}
\def \T{{\bf T}}
\def
\C{{\bf C}}
\def \AA{{\mathcal A}}
\def \BB{{\mathcal B}}
\def
\CC{{\mathcal C}}
\def \DD{{\mathcal D}}
\def \EE{{\mathcal E}}
\def
\FF{{\mathcal F}}
\def \GG{{\mathcal G}}
\def \HH{{\mathcal H}}
\def
\II{{\mathcal I}}
\def \JJ{{\mathcal J}}
\def \KK{{\mathcal K}}
\def
\LL{{\mathcal L}}
\def \MM{{\mathcal M}}
\def \NN{{\mathcal N}}
\def
\OO{{\mathcal O}}
\def \PP{{\mathcal P}}
\def \QQ{{\mathcal Q}}
\def
\RR{{\mathcal R}}
\def \SS{{\mathcal S}}
\def \TT{{\mathcal T}}
\def
\UU{{\mathcal U}}
\def \VV{{\mathcal V}}
\def \ZZ{{\mathcal Z}}
\def
\Pxh{P^{x/h}}
\def \Exh{E^{x/h}}
\def \Px{P^{x}}
\def \Ex{E^{x}}
\def
\Prh{P^{\rho/h}}
\def \Erh{E^{\rho/h}}
\def \p{p_{t}(x,y)}
\def
\({\left(}
\def \){\right)}
\def \lk{\left[}
\def \rk{\right]}
\def
\lc{\left\{}
\def \rc{\right\}}
\def \bsq{\ $\Box$}
\def
\nn{\nonumber}
\def \Bo{\bigotimes}
\def \bo{\times}
\def
\ot{\times}
\def
\bs{\begin{slide} }
\def \es{\end{slide} }
\def \bpr{\begin{proof} }
\def \epr{\end{proof} }
\def \cd{\,\cdot\,}
\def
\st{\stackrel{def}{=}}
\def \as{almost surely }
\def \fix {{\bf !!!!! }}
\def \stl{\stackrel{law}{=}}
\def \std{\stackrel{dist}{=}}
\def \stdto{\stackrel{dist}{\longrightarrow}}
\def  \enc{{ \bf end of current changes}}
\def\square{{\vcenter{\vbox{\hrule height.3pt
                      \hbox{\vrule width.3pt height5pt \kern5pt
                         \vrule width.3pt}
                      \hrule height.3pt}}}}
\def\qed{{\hfill $\square$ }}

             \def \tb{|\!|\!|}

          \eqnsection
\bibliographystyle{amsplain}

\title{  CLT for $L^{ p}$ moduli of continuity  of Gaussian processes
}
          \author{ Michael B. Marcus\,\, Jay Rosen \thanks{Research of
both authors supported by  grants from the National Science
Foundation and PSCCUNY.}}


\maketitle
\eqnsection

\bibliographystyle{amsplain}

            \def \tb{|\!|\!|}

\begin{abstract}   Let  $G=\{G(x),x\in R^1\}$ be a mean zero
Gaussian processes with stationary increments  and set $\si
^2(|x-y|)=  E(G(x)-G(y))^2$. Let  $f$ be a symmetric  function with
$Ef(\eta)<\ff$, where $\eta=N(0,1)$. When $\si^2(s)$ is concave or
when
$\si^2(s)=s^r$, $1<r\leq 3/2$
\[
  \lim_{h\downarrow 0}{ \int_a^bf\(\frac{G(x+h)-G(x)}{\si (h)}
\)\,dx - (b-a)Ef(\eta)\over \sqrt{\Phi(h,\si(h),f,a,b)}}\stl N(0,1)
\] where  $\Phi(h,\si(h),f,a,b)$ is the variance of the numerator.
 This result  continues to hold when $\si^2(s)=s^r$, $3/2<r<2$, for certain
functions $f$, depending on the nature of the  coefficients in their
Hermite polynomial expansion.

The   asymptotic behavior of
$\Phi(h,\si(h),f,a,b)$ at zero, is described in a very large number of cases. 
\end{abstract}

            \section{Introduction}\label{sec-intro}

Let  $G=\{G(x),x\in R^1\}$ be a mean zero  Gaussian  process
with stationary increments, and set
\bea
       E(G(x)-G(y))^2&=&\si ^2( x-y ) \label{m1.2}\\ &=&\si
^2(|x-y|).\nn
\eea
        Clearly $\si^2(0)=0$. To avoid trivialities we assume that
$\si^2(h)\not\equiv 0$.

   When $G$ is continuous and
$\si^2(h)$ is concave for
$h\in[0,h_0]$ for some
$h_0>0$ and satisfies some other very weak conditions,  or when
$\si^2(h)=h^r$, $1<r<2$,   for
$h\in[0,h_0]$, we show in
\cite{lp}   that
\be
\lim_{h\downarrow 0} \int_a^b\bigg|\frac{G(x+h)-G(x)}{\si
(h)}\bigg|^p\,dx =E|\eta |^p (b-a) \label{m1.1}
\ee for  all
$a,b\in R^1$, almost surely, where $\eta$ is a normal random
variable with mean zero and variance one, (sometimes also denoted
by $N(0,1)$).

Obviously, the right-hand side of (\ref{m1.1}) is the expected value
of the integral on the left-hand side for all
$h>0$. Thus one can think of (\ref{m1.1}) as a Strong Law of Large
Numbers for the functional
\be
\int_a^b\bigg|\frac{G(x+h)-G(x)}{\si (h)}\bigg|^p\,dx.\label{m1.3}
\ee
   It is natural to ask  if this functional also satisfies a Central Limit
Theorem  because this would give the next  order in the description
of the  asymptotic behavior of (\ref{m1.3}).

We consider this question in a more general setting.
Fix
$-\ff<a<b<\ff$. Let $d\mu( x)=( 2\pi)^{ -1/2}\exp ( -x^{ 2}/2)\,dx$
denote standard Gaussian measure on
$R^{ 1}$. For any  symmetric  function
$f\in L^{ 2}(R^{ 1},\,d\mu )$, i.e., $Ef(\eta)<\ff$, define
\be I(f,h)=I_G(f,h;a,b)=\int_a^b
f\(\frac{G(x+h)-G(x)}{\si(h)}\)\,dx.\label{m2.5}
\ee
  We obtain CLTs for the functionals
$ I(f,h)$. Clearly they apply to (\ref{m1.3}) by taking
$f(\cd)=|\cd|^p$.

\bt\label{theo-GCLT}
Assume that either $\si^2(h)$ is
concave or that $\si^2(h)=h^r$, $1<r\leq 3/2$. Then for all symmetric
functions $f\in L^{ 2}(R^{ 1},\,d\mu)$
\be
\lim_{h\downarrow 0} {I_G(f,h;a,b)-(b-a)E f(\eta) \over
\sqrt{\mbox{Var
$I_G(f,h;a,b)$}}}
\stackrel{law}{=} N( 0,1).  \label{j.1}
\ee
\et

  When $\si^2(h)=h^r$, $3/2<r<2$ we no longer get (\ref{j.1}) for all symmetric $f\in
L^{ 2}(R^{ 1},\,d\mu)$. However, we do get it for certain $f\in L^{ 2}(R^{
1},\,d\mu)$ depending on the coefficients of the Hermite polynomial
expansion of $f$. 
  Let $\{H_{ m}( x)\}_{m=0}^\ff$ denote the
\label{pageII}Hermite polynomials. (They are
         an orthonormal basis for $L^{ 2}(R^{ 1},\,d\mu )$.) Then for
  symmetric
$f\in L^{ 2}(R^{ 1},\,d\mu )$,
\begin{equation} f( x)=\sum_{ m=0}^{ \ff}a_{2m}H_{2m}(
x)\hspace{ .2in}\mbox{ in   }\,\,L^{ 2}(R^{ 1},\,d\mu ),
\label{a18.1}
\end{equation} where
\be a_{2m}=\int f( x)H_{2m}( x)\,d\mu( x)\label{am.m1}
\ee
       and
\begin{equation}
\sum_{ m=0}^{ \ff}a_{2m}^{ 2}=  \int |f( x)|^{ 2}\,d\mu( x)  <\ff.
\label{a18.2}
\end{equation}

\bt\label{theo-GCLT2} Let  $f\in L^{ 2}(R^{ 1},\,d\mu)$ be symmetric and let
\be k_{ 0}=\inf_{m\ge 1}\{m|a_{2m}\ne 0 \}.\label{bmm.3.1}
\ee
Assume that $\si^2(h)=h^r$,   $0<r\le 2-1/{(2k_0)} $. Then  (\ref{j.1})
 holds.
\et

Clearly, Theorem \ref{theo-GCLT} contains this result when $k_0=1$ but not
when $k_0>1$. We can show that when $f\in L^{ 2}(R^{ 1},\,d\mu)$ is
symmetric and its Hermite polynomial expansion is such that (\ref{bmm.3.1})
holds and   
$\si^2(h)=h^r$, $r>  2-1/{(2k_0)} $,   left-hand side of (\ref{j.1})
converges to a $2k_0$-th order Gaussian chaos. We plan to address this in a
subsequent paper.

\medskip Theorems \ref{theo-GCLT} and \ref{theo-GCLT2} are consequences
of the following   general CLT for $I_G(f,h;a,b)$  and its simple  corllary,
Corollary
\ref{CCLT}.   For
$x,y\in R^{ 1}$ let
\bea
\rho_{h}(x,y)&=&{ 1\over \sigma^{ 2}(h)}E(G(x+h) - G(x))(G(y+h) -
G(y))
\label{ag2.3}\\
      &=&{ 1\over 2\sigma^{
2}(h)}\(\si^2(x-y+h)+ \si^2(x-y-h)-2\si^2(x-y)
  \)\nn\\ &:=&\rho_h(x-y)=\rho_h(y-x).\nn
\eea

\bt\label{BCLT}     Assume
that  for all
$j\in N$
\begin{equation}
\sup_{ a\leq x\leq b}\int_a^b|\rho_h(x-y)|^j\,dy\leq C_{ j}
\int_a^b\!\!\int_a^b|\rho_h(x-y)|^j\,dx\,dy\label{bneed.1aa}
\end{equation} 
where $C_j$ is a constant which can depend on $j$. Assume, furthermore,
that   for all
$j\in N$  
\begin{equation}
\(\int_a^b\!\!\int_a^b|\rho_h(x-y)|^j\,dx\,dy\)^{1/j}
=o\(\int_a^b\!\!\int_a^b
|\rho_h(x-y)|^{j+1}\,dx\,dy\)^{1/(j+1)}.\label{bneed.2r}
\end{equation}
  Then for all symmetric functions   $f\in L^{ 2}(R^{ 1},\,d\mu)$ 
\be
\lim_{h\downarrow 0} {I_G(f,h;a,b)-(b-a)E f(\eta) \over
\sqrt{\mbox{Var
$I_G(f,h;a,b)$}}}
\stackrel{law}{=} N( 0,1).  \label{b3.1www}
\ee
\et

To complete this analysis we need to describe the behavior of  Var
$I_G(f,h;\newline a,b)$  as $h$ decreases to zero. We do this in Sections
\ref{sec-con} and \ref{sec-sigsq}, with varying degrees of precision, depending
on the the function
$\si^2(h)$. We show on page \pageref{1.6gq} that
\begin{equation}
\mbox{Var
$I_G(f,h;a,b)$}=
\sum_{ k=1}^{ \ff}a_{ 2k}^{ 2} \int_{ a}^{ b}\!\!\int_{ a}^{ b}
(\rho_{h}(x-y))^{2 k}
\,dx\,dy .\label{j.5}
\end{equation}
The following table gives the behavior of the integrals in (\ref{j.5}) as $h$
decreases to zero    for many examples of 
$\si^{ 2}( h)$.

 \newpage

\begin{center}
{\bf Table 1}

\vspace{.2in}
\begin{tabular}{ c c| c  }
  &$\si^2(h)$& $\int_{ a}^{ b}\!\!\int_{ a}^{ b}
(\rho_{h}(x-y))^{ k}\,dx\,dy$ \\  \hline   &\\
  1)&$h^r$,\, $r>2-1/k$ &  $\sim C_{ 1,k}\,h^{
(2-r)k}$\\&\\
\hline &\\  2)&$h$ &  $\sim  \frac{2(b-a)}{2k+1}\,h$\\&\\
\hline &\\
  3)&$h^r$,\, $r=2-1/k$, \,$k\ge 2$  &  $\sim C_{
3,k}\,h\,\log 1/h$\\&\\
\hline &\\
4)&  $h^r$,\, $0<r<2-1/k$
&$\sim C_{ 4,k}\,h$\\&\\
\hline &concave\\5)&   regularly varying&$\approx h$\\&strictly positive
index\\
    \hline &\\6)&$\exp(-(\log1/h)^\ga)$, $0<\ga<1$& $\approx \displaystyle
\frac {h}{(\log
1/h)^{k(1-\ga)}}$
\\&\\
    \hline &\\7)&$(\log 1/h)^{-q}$, $q>0$&$\approx
\displaystyle \frac {h}{(\log 1/h)^{k }}$
\end{tabular}
\end{center}

\smallskip \noindent where  
\bea
C_{1,k}&=& {2 r^{ k}| r-1|^{
k}(b-a)^{ (r-2)k+2}\over 2^k((r-2)k+1)((r-2)k+2)}\label{constants}\\
C_{3,k}&=& 2(b-a)\Big|{ r(
r-1)\over 2}\Big|^{k}\nn\\
C_{4,k}&=& 2(b-a)
\int_0^{\ff}
\Big| {|s+1|^{ r} +|s-1|^{ r}-2|s|^{ r}  \over 2} \Big |^{ k}
\,ds\nn.
\eea
\medskip  We use $f\approx g$ at zero, and say that $f$  is
approximately equal to $g$ at zero, to indicate that there exists
constants
$0<C_1\le C_2<\ff$ such that   $C_1\le \liminf_{x\to0}
\frac{f(x)}{g(x)}\le
\limsup_{x\to0}\frac{f(x)}{g(x)}\le C_2$, and  $f\sim g$ at zero,
and say that $f$  is asymptotic
       to $g$ at zero,    to indicate that there exists a constant
$0<C< \ff$ such that  $\lim_{x\to0}\frac{f(x)}{g(x)}= C$. Analogous
definitions apply at infinity. 

\medskip In order to use Table 1 for a given $ f\in L^{ 2}(R^{ 1},\,d\mu)$ it is
necessary to know $k_0$ in (\ref{bmm.3.1}). For the functionals in
(\ref{m1.3}), which were the motivation for this paper, $k_0=1$, since for
these functionals 
$a_2=E(|\eta|^p|\eta^2-1|)/\sqrt2>0$. 
We get the followwing immediate corollary of Theorem \ref{theo-GCLT}.

\bc\label{theo-1} Let $G=\{G(x), x\in [a,(1+\ep)b]\}$, for some $\ep>0$, be a
Gaussian process with stationary increments with increments variance
$\si^2(h)$  that is
concave on  $[0,2(b-a)]$, or satisfies
$\si^2(h)=h^r$,  $1<r\le 3/2$, on  $[0,2(b-a)]$.   Then for all $p\ge 1,$
\be
 \lim_{h\downarrow 0}{ \int_a^b\Big|\frac{G(x+h)-G(x)}{\si
(h)} \Big|^p\,dx -E|\eta |^p (b-a)\over \sqrt{\Phi(h)}}\stl N(0,1)
\label{m1.1qa},
\ee
where $\Phi(h)$ is the variance of the numerator.
\ec

  The follwing table gives the asymptotic behavior of $\Phi(h)$ at
zero for different values of $\si^2(h)$:  
\newpage\begin{center}
{\bf Table 2}

\vspace{.1in}
\begin{tabular}{l  c| c  }   
 &$\si^2(h)$& $\Phi(h)$ \\  \hline   &\\  
(1)& $h^{3/2}$ &  $\sim\({b-a\over \sqrt2}\({3\over
8}\)^2\,\(E( |\eta|^p(\eta^2-1) )\)^{ 2}\)h\log 1/h$\\&\\
\hline &\\(2)& $\stackrel{\mbox{$h^r$}}{ 0<r<3/2}$ &$\stackrel{\mbox{$\sim
2(b-a)\,h\sum_{k=1}^{\ff} \(\(E( |\eta|^pH_{2k}(\eta))\)^{
2}\right.$ }}{\hspace{.8in}
\cd  \left.\int_0^{\ff}
\Big| {|s+1|^{ r} +|s-1|^{ r}-2|s|^{ r}  \over 2} \Big |^{ 2 k}
\,ds\)}$\\&\\\hline &\\(3)&$h$&$\sim2(b-a)\,h\sum_{k=1}^{\ff}\(E(
|\eta|^pH_{2k}(\eta))\)^{ 2}\, \frac{1}{2k+1}$\\&\\
\hline & concave& \\(4)&regularly varying&$\approx h$\\&strictly  positive
index\\
 \hline &  &\\ (5)&concave slowly varying&$\approx 
\displaystyle\({h\si'(h)\over
\si(h )}\)^2 \,h$  
\end{tabular}

\vspace{.1in} 
\end{center}
\noindent  with the additional  condition, in the final expression, that
$h{d\over dh}\si^2(h)$ is increasing.   It is easy to see that the
last  entry in this table agrees with   (6) and (7), with $k=2$, in Table 1.

In \cite[Theorem 2.2]{ST} Sodin 
and Tsirelson give a general CLT for Gaussian functionals which gives
some, but not all, of the cases covered by Theorem \ref{BCLT}.   
Their theorem states that (\ref{b3.1www}) holds whenever  
       \be
\lim_{h \to 0}\sup_{a\leq x\leq b }\int_{ a}^{
b}|\rho_{h}(x-y)|\,dy=0\label{am1.6}
\ee and
   for all
$k\in N$
\begin{equation}
\liminf_{h\downarrow 0 }{ \int_{ a}^{ b}\int_{ a}^{
b}|\rho_{h}(x-y)|^{2 k}\,dx\,dy\over
\sup_{a\leq x\leq b }
\int_{ a}^{ b}|\rho_{h}(x-y)|\,dy}>0.\label{a3.2}
\end{equation}
  When $f$ is increasing on $[0,\ff)$ it suffices to have
(\ref{a3.2}) for $k=1$.

For all  the examples in Table 1 we have that for all
$k\in N$
\begin{equation}
\sup_{a\leq x\leq b }\int_{ a}^{
b}|\rho_{h}(x-y)|^{k}\,dy\approx  \int_{ a}^{ b}\!\!\int_{ a}^{ b}
|\rho_{h}(x-y)|^{ k}\,dx\,dy\label{j.6}
\end{equation}
so that (\ref{am1.6}) holds for all these examples and
condition (\ref{a3.2}) for
$k=1$ is equivalent to
\begin{equation}
\liminf_{h\downarrow 0 }{ \int_{ a}^{ b}\int_{ a}^{
b}|\rho_{h}(x-y)|^{2 }\,dx\,dy\over
  \int_{ a}^{ b}\int_{ a}^{
b}|\rho_{h}(x-y)|\,dx\,dy}>0.\label{a3.2j}
\end{equation}
It is easily seen that this holds in case $5)$ of Table 1 but not in cases
$6)$ and $7)$, nor  when $\si^2(h)=h^{ r}$ for $1<r\leq 3/2$.    Actually,
the CLT in \cite[Theorem 2.2]{ST}, as  it applies to $I_G(f,h;a,b)$, is
contained in Corollary \ref{CCLT} with $k_0=2$.

  It should be clear that we can not get
classical CLTs for
 $I_G(f,h;a,b)$  for all Gaussian process. For example      when
$\si^2(h)=h^2$,
\be
       \int_{ a}^{ b}\int_{ a}^{ b} |\rho_{h}(x-y)|^{ 2k}
\,dx\,dy=2^k(b-a)^2,
\ee so that  $ \lim_{h\downarrow 0}\mbox{Var $I_G(f,h;a,b)$}\ne
0$.  To make this example more explicit suppose that  Gaussian process
$G$ is integrated Brownian motion, then 
\be
\lim_{h\downarrow 0} \int_a^b\bigg|\frac{G(x+h)-G(x)}{h}\bigg|^p\,dx
=\int_a^b | B(x)|^p\,dx\qquad \mbox{a. s.}
\label{m1.1ac}
\ee where $B$ is Brownian motion.   
Obviously, the right-hand side of (\ref{m1.1ac}) is not $N(0,1)$.

  \medskip    In Section \ref{sec-4} we prove the general
Theorem \ref{BCLT}  and Corollary \ref{CCLT}. To obtain Theorems
\ref{theo-GCLT} and 
\ref{theo-GCLT2} we must verify that the conditions of Theorem \ref{BCLT}
and  Corollary \ref{CCLT} hold, when
$\si^2(h)$ is concave or when $\si^2(h)=h^{ r}$ for   $1<r\leq 2-1/(2k_0)$.
In  Section \ref{sec-con} we 
do this for
$\si^2(h)$  concave
and in Section \ref{sec-sigsq}  when $\si^2(h)$   is a
power.  We give the proofs of Theorems \ref{theo-GCLT} and
\ref{theo-GCLT2} in Section
\ref{sec-proofs} and also point out how we obtain the estimates in Tables 1
and 2.

\section{Proof of Theorem \ref{BCLT}}\label{sec-4}

      Let $\phi_h(x-y)=\si^2(h)\rho_h(x-y)$.   Note that
\be
\phi_h(x-y) ={1 \over 2}\(  \si^2(x-y+h)+\si^2(x-y-h)-2\si^2(x-y)\)
\ee
       The $2k$-th Wick product for a mean zero Gaussian random
variable
$Z$ is
\begin{equation}
        :Z^{ 2k}:\,=\sum_{ j=0}^{ k}( -1)^{ j}{ 2k\choose 2j}E( Z^{
2j})\,\,Z^{ 2( k-j)}.\label{7.50}
\end{equation}
If $Z=N( 0,1)$ then  $:Z^{
2k}:=\sqrt{(2k)!}H_{2k}( Z)$.    Hence if
$\si^{ 2}_{ Z}$ denote the variance of
$Z$,
\be
        :\({Z \over \si_{ Z}}\)^{ 2k}:=\sqrt{(2k)!}H_{2k}\({Z \over \si_{
Z}}\).
\label{7.54a}
\ee

\bl\label{theo-weak2m} Let $G$ be a mean zero Gaussian process
with stationary increments.    Assume that
\begin{equation}
\sup_{ a\leq x\leq b}\int_a^b|\rho_h(x-y)|^{2k}\,dy\leq C
\int_a^b\!\!\int_a^b|\rho_h(x-y)|^{2k}\,dx\label{mneed.1aa}
\end{equation} and for all $j<2k$
\begin{equation}
\sup_{ a\leq x\leq b}\(\int_a^b|\rho_h(x-y)|^j\,dy\)^{1/j}
=o\(\int_a^b\!\!\int_a^b
|\rho_h(x-y)|^{2k}\,dx\,dy\)^{1/2k}.\label{mneed.2r}
\end{equation}
Then
\begin{equation} {\lim_{ h\downarrow 0}  \int_{ a}^{ b} :\(\frac{G(
x+h)- G( x) }{\si ( h)}\)^{2k}:\,dx
\over
\sqrt{\int_{ a}^{ b}\int_{ a}^{ b} |\rho_{h}(x,y)|^{ 2k}
\,dx\,dy}}\stackrel{law}{=} \sqrt{(2k)!}\,\,N( 0,1). \label{m3.1}
\end{equation}
\el

\Proof We write
\be
       {\int_{ a}^{ b} :\(\frac{G( x+h)- G( x) }{\si ( h)}\)^{2k}:\,dx
\over
\sqrt{\int_{ a}^{ b}\int_{ a}^{ b} |\rho_{h}(x,y)|^{ 2k}
\,dx\,dy}}=  {\int_{ a}^{ b} :(G( x+h)- G( x))^{ 2k}:\,dx \over
\(\int_a^b\!\!\int_a^b |\phi_h(x-y)|^{ 2k}\,dx\,dy\)}
\ee and  show that for each $n\ge 1$
\begin{eqnarray} &&
\lim_{ h\rar 0}E\(\lc
       { \int_{ a}^{ b} :(G( x+h)- G( x))^{ 2k}:\,dx \over
\(\int_a^b\!\!\int_a^b |\phi_h(x-y)|^{ 2k}\,dx\,dy\)^{1/2}}\rc^{
n}\)\nn\\ &&\hspace{ 1.5in} =\left\{\begin{array}{ll} {( 2m)! \over
2^{ m}m!}\((2k)!\)^{m} &\mbox{ if }n=2m\\ 0&\mbox{ otherwise.}
\end{array}
\right.
\label{7.53t}
\end{eqnarray} Since the right-hand side of (\ref{7.53t}) are the
moments    of the  right-hand side of (\ref{m3.1}) the theorem is
proved.

It follows from \cite[Lemma 2.2]{mem} that
\begin{equation}  E\(\prod_{ i=1}^{ n}:(G( x_{ i}+h)- G( x_{ i}))^{
2k}:\)=
\sum_{\pi\in \mathcal{P}}\(\prod_{ (i,i')\in \pi}\phi_{ h}( x_{
i}-x_{ i'})\)\label{7.54mem}
\end{equation} where the sum runs over all pairings
$\pi\in \mathcal{P}$, the set of pairings of the $2kn$ elements
which consist of $2k$ copies of each of the letters $x_{ i},\,1\leq
i\leq n$, subject to the   restriction that  no single letter
$x_{ i}$ is  paired with itself.

      We  say that the letters  $x_{ i},\,x_{ j}$ are connected in the
pairing
$\pi$ if we can find some sequence $( i_{ m},i_{
m+1}),\,m=1,\ldots$ of pairs in
$\pi$ with $ i_{ 1}=i,i_{p}=j$ for some $p$. By decomposing the set
of letters $x_{ i},\,1\leq i\leq n$ into connected components we can
write (\ref{7.54mem}) as
\begin{eqnarray} && E\(\prod_{ i=1}^{ n}:(G( x_{ i}+h)- G( x_{
i}))^{ 2k}:\)\label{7.54}\\ &&\qquad   = \sum_{ l=1}^{ [n/2]}
\sum_{ C_{ 1}\cup C_{ 2}\cup\cdots \cup C_{ l} =\{\,x_{ i},
i=1,\ldots,n \}}\prod_{ j=1}^{ l}
\sum_{\pi\in \mathcal{P}( C_{j })}\(\prod_{ (i,i')\in \pi}\phi_{ h}(
x_{ i}-x_{ i'})\)\nonumber
\end{eqnarray}  where the second sum runs over all partitions of
$\{\,x_{ i}, i=1,\ldots,n \}$  into $l$ sets, $C_1\ldots,C_l$ with
       $|C_{ i}|\geq 2$, $i=1,\ldots,l$. ($|C|:= \#  \mbox{ of elements
in  } C$.) The third sum runs over all pairings
$\pi\in \mathcal{P}( C_{j })$, the set of pairings of the set of
$2k|C_{j }|$ elements which consists of $2k$ copies of each of the
letters $x_{ i}\in C_{j }$, subject to the following two restrictions:
\begin{itemize}
\item[(i)]  no single letter
$x_{ i}$ is  paired with itself;
\item[(ii)]  for any partition $C_{j }=A\cup B$, at least one letter of
$A$, is paired with a letter of $B$. ( This condition states that $C_j$ can
not be further decomposed into  connected components.)
\end{itemize}

      We show below that the only  non-zero terms  of the left-hand
side of (\ref{7.53t}) comes when 
  $n=2m$  and the  partitions have
       $m$ parts,
$(C_{ 1},C_{ 2},\cdots, C_{ m} )$, in which case   all parts
necessarily have
      two elements; that is,   from pairings of $\{\,x_{ i}, i=1,\ldots,2m
\}$.
   Referring again to   (\ref{7.54mem})
we see that for each partition of this sort
\be
\prod_{ (i,i')\in \pi}\phi_{ h}( x_{
i}-x_{ i'})=\prod_{j=1}^n \phi_{ h}^{2k}( x_{
i_j}-x_{ i'_{j'}})
\ee
where $C_j=(i_j,i'_{j'})$.

       Since there are ${( 2m)! \over 2^{ m}m!}$ pairings of
$\{\,x_{ i}, i=1,\ldots,2m \}$  and $(2k!)$ ways to arrange the two sets of
$2k$  elements in each pairing, it follows from (\ref{7.54})  that
\begin{eqnarray} \lefteqn{ E\(\lc
       \int_{ a}^{ b} :(G( x+h)- G( x))^{ 2k}:\,dx \rc^{
2m}\)\label{7.57}}\\&& ={( 2m)! \over 2^{ m}m!}\((2k)!\)^{m}
\(\int_a^b\!\!\int_a^b |\phi_h(x-y)|^{ 2k}\,dx\,dy\)^{m}
\nn\\ &&\qquad +\sum_{ l=1}^{ m-1}
\sum_{ C_{ 1}\cup C_{ 2}\cup\cdots \cup C_{ l} =\{\,x_{ i},
i=1,\ldots,2m
\}}\nn\\ &&\hspace{1in}\int_{ [a,b]^{ 2m}}\prod_{ j=1}^{ l}
\sum_{\pi\in \mathcal{P}( C_{j })}\(\prod_{ (i,i')\in \pi}\phi_{ h}(
x_{ i}-x_{ i'})\)\,\prod_{i=1}^{2m}dx_i.\nonumber
\end{eqnarray}

Since the first term to the right of the equal  sign in (\ref{7.57})
gives (\ref{7.53t}), and
\bea &&\int_{ [a,b]^{ 2m}}\prod_{ j=1}^{ l}
\sum_{\pi\in \mathcal{P}( C_{j })}\(\prod_{ (i,i')\in \pi}\phi_{ h}(
x_{ i}-x_{ i'})\)\,\prod_{i=1}^{2m}dx_i\\ &&\qquad =\prod_{
j=1}^{ l}\int_{ [a,b]^{ |C_j|}}\sum_{\pi\in \mathcal{P}( C_{j })}
\(\prod_{ (i,i')\in \pi}\phi_{ h}( x_{ i}-x_{ i'})\)\,\prod_{x_i\in C_j}
dx_i\nn
\\ &&\qquad =\prod_{ j=1}^{ l}\sum_{\pi\in \mathcal{P}( C_{j
})}\int_{ [a,b]^{ |C_j|}}
\(\prod_{ (i,i')\in \pi}\phi_{ h}( x_{ i}-x_{ i'})\)\,\prod_{x_i\in C_j}
dx_i\nn,
\eea to complete the proof of (\ref{7.53t}), when $n$ is even, it
      suffices to show that for any  set,  say 
$C_{p}$, with
$|C_{p}|\geq 3$, and any $\pi\in \mathcal{P}( C_p)$
\begin{eqnarray} &&
\int_{ [a,b]^{ |C_p|}}\prod_{ (i,i')\in \pi}
\phi_{ h}( x_{i}-x_{ i'})\,\prod_{x_i\in C_p}dx_i\label{7.58}\\
&&\qquad =o\(\(\int_a^b\!\!\int_a^b |\phi_h(x-y)|^{
2k}\,dx\,dy\)^{|C_p|/2} \).\nn
\end{eqnarray}

      To   obtain (\ref{7.58})  choose any pair of letters $x_{i},x_{ i'}$
with
$(i,i')\in
\pi$. Suppose that 
$j$ is the number of times that $(i,i')$ occurs in $\pi$, then we must
have
$1\leq j< 2k$, since  if $j=2k$   restriction (ii) would be violated.
Each variable
$x_{ r}$ on the   left-hand side of (\ref{7.58}) occurs  precisely $2k$
times.
       Pick such an
$x_r\neq x_i $ or $x_{i'}$ and use the generalized H\"older's
inequality together with (\ref{mneed.1aa}) to obtain the bound
\begin{eqnarray} &&
\sup_{ a\leq d_{ j}\leq b,\,\forall j}\int_{a}^{b}\prod_{j=1}^{ 2k}
\phi_{ h}(x_{r}-d_{ j})\,dx_{r}\label{7.58j}\\ &&\qquad
\leq \sup_{ a\leq d_{ j}\leq b,\,\forall j}\prod_{j=1}^{
2k}\(\int_{a}^{b}|
\phi_{ h}(x-d_{ j})|^{ 2k}\,dx\)^{ 1/2k}\nonumber\\ && \qquad
\leq\prod_{j=1}^{ 2k} \sup_{ a\leq d_{ j}\leq b}\(\int_{a}^{b}|
\phi_{ h}(x-d_{ j})|^{ 2k}\,dx\)^{ 1/2k}\nonumber\\&&\qquad
\leq \sup_{ a\leq d\leq b}\(\int_{a}^{b}|
\phi_{ h}(x-d)|^{ 2k}\,dx\)\nonumber\\ &&\qquad
\leq K\int_a^b\!\!\int_a^b |\phi_h(x-y)|^{ 2k}\,dx\,dy.\nonumber
\end{eqnarray}
      Here $\{d_j\}_{j=1}^{2k}$ represents the different elements
$x_j\in C_p$ that
$x_r$ is paired with. Several of the $d_j$ may be the same.

We proceed to  successively  bound the integrals over each
$x_{ p}$ with
$p\neq i,i'$. Now, however, there may be less than $2k$ remaining
factors containing $x_{ p}$ since some factors may have been
bounded at an earlier stage. If, say there are $q$ factors left when
we  bound $x_{ p}$, then as in (\ref{7.58j}) we obtain
\begin{eqnarray} &&
\sup_{ a\leq d_{ j}\leq b,\,\forall j}\int_{a}^{b}\prod_{j=1}^{q}
\phi_{ h}( x_{p}-d_{ j})\,dx_{p}\label{7.58k}\\ &&\qquad 
\leq (b-a)^{1-q/2k}\sup_{ a\leq d_{ j}\leq b,\,\forall
j}\prod_{j=1}^{ q}\(\int_{a}^{b}|
\phi_{ h}( x_{p}-d_{ j})|^{ 2k}\,dx_{ p}\)^{ 1/2k}\nonumber\\
&&\qquad
\leq (b-a)^{1-q/2k}\sup_{ a\leq d\leq b}\(\int_{a}^{b}|
\phi_{ h}( x_{p}-d)|^{ 2k}\,dx_{ p}\)^{ q/2k}\nonumber\\
&&\qquad
\leq(b-a)^{1-q/2k}K\(\int_a^b\!\!\int_a^b |\phi_h(x-y)|^{
2k}\,dx\,dy\)^{ q/2k}.\nonumber
\end{eqnarray}

       Note that the number of pairs in any $\pi\in \mathcal{P}( C_p)$ is
$|C_p|k$. Thus we see that after bounding successively all the integrals
involving
$x_{ r}$ with
$r\neq i,i'$ we have  for some $1\le j<2k$
\begin{eqnarray} \lefteqn{
\int_{ [a,b]^{ |C_p|}}\prod_{ (i,i')\in \pi}
\phi_{ h}( x_{i}-x_{ i'})\,\prod_{x_i\in C_p}dx_i\label{7.58m}}\\
&&\!\leq K'\(\int_a^b\!\!\int_a^b |\phi_h(x-y)|^{
2k}\,dx\,dy\)^{(|C_p|k-j)/2k}
\(\int_a^b\!\!\int_a^b |\phi_h(x_{ i}-x_{ i'})|^{ j}\,dx_{ i}\,dx_{
i'}\)\nn 
\end{eqnarray}
      where $K'<\ff$ does not depend on $h$. Since by
(\ref{mneed.2r})
\be
\int_a^b\!\!\int_a^b |\phi_h(x_{ i}-x_{ i'})|^{ j}\,dx_{ i}\,dx_{
i'}=o\(\int_a^b\!\!\int_a^b |\phi_h(x_{ i}-x_{ i'})|^{ 2k}\,dx_{
i}\,dx_{ i'}\)^{j/2k},
\ee we get (\ref{7.58}).

 At this point it should be clear that (\ref{7.53t}) is zero when $n$ is odd
since any partition of $\{\,x_{ i}, i=1,\ldots,n \}$  into $l$ sets, $C_1\ldots,C_l$
with
       $|C_{ i}|\geq 2$, $i=1,\ldots,l$,  contains at least one set with three
of more elements.
\qed

\medskip  To proceed we  need some more information
about the Hermite polynomial expansion of functions in $ L^{ 2}(R^{
1},\,d\mu )$. It is clear  that
\begin{equation} E( f( X))=\int f( x)\,d\mu( x)=a_{ 0}\label{18.4}
\end{equation} so that
\begin{equation} f(X)-E( f( X))=\sum_{ m=1}^{ \ff}a_{ 2m}H_{2 m}(
X)\hspace{ .2in}\mbox{ in   }\,\,L^{ 2}(R^{ 1},\,d\mu ).
\label{18.5}
\end{equation} Let  $X$ and $Y$  be $N(0,1)$ and let $(X,Y)$ be a
two dimensional Gaussian random variable. Then
\begin{equation} E(H_{ 2m}( X)H_{ 2n}( Y))=(E( XY))^{ 2m}\de_{
m,n}.\label{18.3}
\end{equation} This follows by setting $Y=\al X+(1-\al^2)^{1/2}Z$,
where
$\al=E(XY)$ and
$Z$ is
$N(0,1)$ and is independent of $X$, and using  the  relationship
\begin{equation}
\sum_{ m=0}^{ \ff}{ \la^{ m}\over \sqrt{m!}}H_{ m}( x) =\exp(\la
x-\la^{ 2}/2).\label{18.3g}
\end{equation}
      Consequently, it follows from (\ref{18.3}) that
\begin{equation}
\mbox{ Cov }(f(X),\,f(Y) )=\sum_{ m=1}^{ \ff}a_{ 2m}^{ 2}(E( XY))^{
2m}.\label{18.6q}
\end{equation}

  For each
$h$ we consider the symmetric positive definite kernel
$\rho_{h}(x,y)=\rho_{h}(x-y)$. Note that by stationarity and the
Cauchy--Schwarz inequality
\begin{equation} |\rho_{h}(x-y)|\leq 1\qquad \hspace{
.2in}\forall\,x,y\in R^{ 1}.\label{ag2.3t}
\end{equation}
         Therefore,  by (\ref{18.6q})
\bea \lefteqn{
\mbox{ Var }\(\int_{a}^{ b} f\(\frac{G(x+h)-G(x)}{\si(h)}\)\,dx\)
\label{1.6gq}}\\ && =\int_{a}^{ b}\!\!\int_{a}^{ b}\mbox{ Cov
}\(f\(\frac{G(x+h)-G(x)}{\si(h)}\),\,f\(\frac{G(y+h)-G(y)}{\si(h)}\)
\)\,dx\,dy\nn\\ &&=
\sum_{ m=1}^{ \ff}a_{ 2m}^{ 2} \int_{ a}^{ b}\!\!\int_{ a}^{ b}
(\rho_{h}(x-y))^{2 m}
\,dx\,dy . \nn
\eea

\medskip\noindent{\bf Proof of Theorem \ref{BCLT} }    Clearly   we need
only consider  $f\in L^{ 2}(R^{ 1},\,d\mu)$ of the form
\begin{equation} f( x)=\sum_{ m=k_0}^{ \ff}a_{2m}H_{2m}(
x).
\label{m.6}
\end{equation}  
 To begin suppose
that there are only a finite number of terms in (\ref{m.6}) so that for some
$k_{1}<\ff$
\begin{equation} f ( x)=\sum_{ m=k_0}^{ {k_1}}a_{2m}H_{2m}(
x).
\label{m.6a}
\end{equation}

 Let $Y_{ h}=\int_{ a}^{ b}f \(\frac{G( x+h)- G( x) }{\si ( h)}\)\,dx$.
By (\ref{1.6gq}) we see that
\begin{equation}
\mbox{ Var }(Y_{ h})=
\sum_{ m=k_{ 0}}^{k_1}a_{2m}^{ 2} \int_{ a}^{ b}\!\!\int_{ a}^{
b} (\rho_{h}(x-y))^{2m}
\,dx\,dy .\label{b.6gq}
\end{equation} Since $|\rho_{h}(x-y)|\le 1$ we see that 
\bea 
&&a_{2k_0}^{ 2} \int_{ a}^{ b}\!\!\int_{ a}^{
b} (\rho_{h}(x-y))^{2k_0}\,dx\,dy\label{mend}\\
&&\qquad \quad\le  \mbox{ Var }(Y_{ h})\le \sum_{ m=k_{
0}}^{k_1}a_{2m}^{ 2} \int_{ a}^{ b}\!\!\int_{ a}^{ b}
(\rho_{h}(x-y))^{2k_0}\,dx\,dy \nn.
\eea

We obtain (\ref{b3.1www}) by showing that, in the limit, as $h\downarrow0$, 
the moments of the left-hand side are equal to the moments of the
right-hand side, (as in the proof of Lemma \ref{theo-weak2m}).
 We have
\begin{eqnarray} &&
E\lc\(\sum_{ m=k_{ 0}}^{k_{ 1}}a_{ 2m}\int_{ a}^{ b}H_{2 m}
\(\frac{G( x+h)- G( x) }{\si ( h)}\)\,dx\)^{ n}\rc\label{j.10}\\ &&\qquad 
=\sum_{ m_{ i}=k_{ 0};\,i=1,\ldots,n}^{k_{ 1}}
\(\prod_{i=1}^{ n}a_{ 2m_{ i}}\)\nonumber\\ && 
\hspace{ 1in}\int_{ [a,b]^{ n}} E\lc \prod_{i=1}^{ n}
H_{2 m_{ i}}\(\frac{G( x_{ i}+h)- G( x_{ i}) }{\si ( h)}\)\rc \prod_{i=1}^{
n}\,dx_{ i}.
\nonumber
\end{eqnarray}
  As in the proof of Theorem \ref{theo-weak2m} we have
\begin{eqnarray} &&
  E\lc \prod_{i=1}^{ n}
H_{2 m_{ i}}\(\frac{G( x_{ i}+h)- G( x_{ i}) }{\si ( h)}\)\rc \label{j.11}\\
&&\qquad  =\sum_{ l=1}^{ [n/2]}
\sum_{ C_{ 1}\cup C_{ 2}\cup\cdots \cup C_{ l} =\{\,x_{ i},
i=1,\ldots,n \}}\prod_{ j=1}^{ l}
\sum_{\pi\in \mathcal{P}( C_{j })}\(\prod_{ (i,i')\in \pi}\rho_{ h}(
x_{ i}-x_{ i'})\)\nonumber
\end{eqnarray}  where the second sum runs over all partitions of
$\{\,x_{ i}, i=1,\ldots,n \}$  into $l$ sets, $C_1\ldots,C_l$ with
       $|C_{ i}|\geq 2$, $i=1,\ldots,l$. ($|C|:= \#  \mbox{ of elements
in  } C$)
and if $C=\{ x_{ 1}, \dots, x_{k }\}$,
then $\mathcal{P}( C)$ is the set of pairings of the $\sum_{ i=1}^{
k}2m_{ i}$ elements consisting of
$2m_{ i}$ copies of the letter $x_{ i}$ subject to the same two
retrictions as in the proof of Lemma \ref{theo-weak2m}:
\begin{itemize}
\item[(i)]  no single letter
$x_{ i}$ is  paired with itself;
\item[(ii)]  for any partition $C=A\cup B$, at least one letter of $A$,
is paired with a letter of $B$.
\end{itemize} Of course all $k_{ 0}\leq m_{ i}\leq k_{ 1}$.

Let
\begin{equation}
\mathcal{G}=\{ C_{ 1}\cup C_{ 2}\cup\cdots \cup C_{ l} =\{\,x_{ i},
i=1,\ldots,n \}\,|\,|C_{i}| =2,\,i=1,\ldots,l.\}\label{j.12}
\end{equation}
Then necessarily for partitions in $\mathcal{G}$, $n$ is even, $l=n/2$
and the restrictions on $\mathcal{P}( C_{ i})$ show that if $C_{ i}=\{
x_{ i},x_{ j} \}$ then $m_{ 2j}=m_{ 2i}$.  In this case the contribution   to 
the last line of (\ref{j.10}) is 
\begin{equation}
\prod_{ i=1}^{ n/2} \(\int_{ a}^{ b}\!\!\int_{ a}^{b}
(\rho_{h}(x-y))^{2m_{ i}}
\,dx\,dy\) \label{j.13}.
\end{equation}
 There are ${( 2l)! \over 2^{ l}l!}$ pairings of
$\{\,x_{ i}, i=1,\ldots,n=2l \}$.
Hence the contribution of all the
partitions in $\mathcal{G}$ to (\ref{j.10}) is 
\bea
&&
{( 2l)! \over 2^{ l}l!}\sum_{ m_{ i}=k_{ 0};\,i=1,\ldots,n/2}^{k_{ 1}}
\(\prod_{i=1}^{ n/2}a^{ 2}_{ 2m_{ i}} \(\int_{ a}^{ b}\!\!\int_{ a}^{b}
(\rho_{h}(x-y))^{2m_{ i}}
\,dx\,dy\)    \)\label{j.14}\\
&&\qquad ={( 2l)! \over 2^{ l}l!}\prod_{i=1}^{ n/2}\(\sum_{ m_{ i}=k_{
0}}^{k_{ 1}}a^{ 2}_{ 2m_{ i}} \(\int_{ a}^{ b}\!\!\int_{ a}^{b}
(\rho_{h}(x-y))^{2m_{ i}}
\,dx\,dy\)    \)\nn\\
&&\qquad ={( 2l)! \over 2^{ l}l!}\(\mbox{ Var }(Y_{ h})\)^{ l}\nn
\eea
where the last line comes from (\ref{b.6gq}).

  Thus, as in the proof of Lemma \ref{theo-weak2m}, it
suffices to show that
  for any  set  say
$C_{p}$, with
$|C_{p}|\geq 3$, and any $\pi\in \mathcal{P}( C_p)$
\be  
\int_{ [a,b]^{ |C_p|}}\prod_{ (i,i')\in \pi}
\rho_{ h}( x_{i}-x_{ i'})\,\prod_{x_i\in C_p}dx_i =o\(\(\mbox{ Var
}(Y_{ h})\)^{|C_p|/2} \).\label{b7.58}
\ee
 Suppose that $|C_p|=k$. We relabel the elements of $C_p$, 
$x_1,\ldots,x_k$ and choose them so that 
 $m_{1}\leq m_{2}\leq\ldots\leq m_{k}$.
  If   there are no strict inequalities, i.e., if $m_{1}= m_{2}=\ldots=
m_{k}$,  then,  because of (\ref{bneed.2r}), we are in the same situation as 
in the proof of Theorem
\ref{theo-weak2m}  and we obtain  
\be  
\int_{ [a,b]^{ |C_p|}}\prod_{ (i,i')\in \pi}
\rho_{ h}( x_{i}-x_{ i'})\,\prod_{x_i\in C_p}dx_i =o\(\int_{ a}^{ b}\!\!\int_{ a}^{b}
(\rho_{h}(x-y))^{2m_{ 1}}
\,dx\,dy)\)^{|C_p|/2} .\label{b7.58q}
\ee
Using (\ref{mend}) and the fact that $|\rho_h(\cd)|\le 1$, we see that this
implies (\ref{b7.58}).

If there is at least one strict inequality, that is, if $m_{j}<
m_{j+1}$, for at least one $1\leq j\leq k$, it follows from  the second 
restriction  on
$\pi$, that we can find some
$( j,j')\in \pi$ with $m_{j}< m_{j'}$.  Set
\begin{equation}
\| \rho_{h} \|_{ 2m_i}= \(\int_{ a}^{ b}\!\!\int_{ a}^{ b}
(\rho_{h}(x-y))^{2m_i}
\,dx\,dy\)^{ 1/2m_i}.\label{b7.58aa}
\end{equation} Using   (\ref{mend})  again and the fact that $\| \rho_{h}
\|_{ 2m_{i}}^{2 m_{i}}\leq \| \rho_{h} \|_{
2k_0}^{2 k_0}$,  we  see that  to obtain (\ref{b7.58}), it suffices to show that
\begin{eqnarray} &&
\int_{ [a,b]^{ |C_p|}}\prod_{ (i,i')\in \pi}
\rho_{ h}( x_{i}-x_{ i'})\,\prod_{x_i\in C_p}dx_i =o\(\prod_{ i=1}^{
k}\| \rho_{h} \|_{ 2m_{i}}^{ m_{i}} \).\label{b7.58a}
\end{eqnarray}

To show that (\ref{b7.58a}) holds,   we  successively bound the
integrals over $x_{1}, x_{2},\newline\ldots, x_{k}$ using H\"{o}lder's
inequality, as described in the proof of  Theorem
\ref{theo-weak2m}.
This shows that each factor of the form $\rho_{ h}( x_{i}-x_{ i'})$
  with $i< i'$ makes a contribution which is  $O\(\|
\rho_{h}
\|_{ 2m_{i}}\)$. When $m_{i}=m_{i'}$ we can write this as $O\(\|
\rho_{h}
\|_{ 2m_{i}}^{ 1/2}\|
\rho_{h}
\|_{ 2m_{i'}}^{ 1/2}\)$. When  $m_{i}<m_{i'}$ it follows from   
(\ref{bneed.2r}) that   the bound $O\(\|\rho_{h}
\|_{ 2m_{i}}\)=o\(\|\rho_{h}
\|_{ 2m_{i}}^{ 1/2}\|\rho_{h}\|_{ 2m_{i'}}^{ 1/2}\)$. 

Since there are
$2m_{i}$ factors containing
$x_{ i}$ for each $i$, we always get a bound which is
$O\(\prod_{ i=1}^{
k}\| \rho_{h} \|_{ 2m_{i}}^{ m_{i}}\)$. The desired estimate
(\ref{b7.58a})  follows  because, as we pointed out above, for some
$( j,j')\in \pi$,   $m_{j}< m_{j'}$. Thus we get (\ref{b3.1www}) when the
Hermite polynomial expansion of $f$ contains a finite
number of terms.

To remove this restriction consider an $f$ as in (\ref{m.6})   and let
\begin{equation} f_{n}( x)=\sum_{ m=k_0}^{n}a_{2m}H_{2m}(
x) 
\label{m.6r}
\end{equation}
Set $Y_{ 
h}=\int_{ a}^{ b} f  \(\frac{G( x+h)- G( x) }{\si ( h)}\)\,dx$ and $Y_{n, 
h}=\int_{ a}^{ b} f_{n}  \(\frac{G( x+h)- G( x) }{\si ( h)}\)\,dx$. Using
(\ref{1.6gq})  and  the fact that $\int_{ a}^{ b}\!\!\int_{ a}^{ b}
(\rho_{h}(x-y))^{2m}\,dx\,dy $ is
  decreasing as $m$ increases, we have  
\bea && \lim_{ n\rar\ff}\sup_{ h}E\lc\({Y_{ h} -Y_{n, h} \over
\sqrt{\mbox{Var }( Y_{ h})}}\)^{ 2}\rc\label{d3.2tt}\\
&&\qquad  =\lim_{ n\rar\ff}\sup_{ h}{\sum_{
m=n+1}^{\ff}a_{2m}^{ 2} \int_{ a}^{ b}\!\!\int_{ a}^{ b}
(\rho_{h}(x-y))^{2m}
\,dx\,dy  \over \sum_{ m=k_0}^{\ff}a_{2m}^{ 2} \int_{ a}^{
b}\!\!\int_{ a}^{ b} (\rho_{h}(x-y))^{2m}
\,dx\,dy }\nn\\ &&\qquad 
\leq\lim_{ n\rar\ff}\sup_{ h} {\sum_{ m=n+1}^{\ff}a_{2m}^{ 2} \int_{
a}^{ b}\!\!\int_{ a}^{ b} (\rho_{h}(x-y))^{2m}
\,dx\,dy  \over a^{ 2}_{ 2k_{ 0}} \int_{ a}^{ b}\!\!\int_{ a}^{ b}
(\rho_{h}(x-y))^{2k_{ 0}}
\,dx\,dy }\nn\\ &&\qquad \leq \lim_{ n\rar\ff}{ 1\over a^{ 2}_{ 2k_{
0}}}\,\sum_{ m=n+1}^{\ff}a_{2m}^{ 2}= 0. \nn
\eea
Therefore, we can take the weak limit of
\be
\lim_{h\downarrow 0} {I_G(f_{n},h;a,b)-(b-a)E f_{n}(\eta) \over
\sqrt{\mbox{Var
$I_G(f_{n},h;a,b)$}}}
\ee  
as $n  \to\ff$ and obtain (\ref{b3.1www}).\qed

\medskip  We get the following simple corolary of Theorem \ref{CCLT} in
which gives a weaker condition than (\ref{bneed.2r}) when an additional
regularity condition is satisfied. 

\bc\label{CCLT}  Let  $f\in L^{ 2}(R^{ 1},\,d\mu)$ be symmetric and suppose 
that    its Hermite polynomial expansion is such that (\ref{bmm.3.1})  
holds. Assume that  (\ref{bneed.1aa}) holds for all
$  j\in N$.
  Assume, furthermore,
that   for all
$1\le j<2k_0$ 
\begin{equation}
\(\int_a^b\!\!\int_a^b|\rho_h(x-y)|^j\,dx\,dy\)^{1/j}
=o\(\int_a^b\!\!\int_a^b
|\rho_h(x-y)|^{2k_0}\,dx\,dy\)^{1/(2k_0)}\label{bneeq}
\end{equation}
and
\be
\liminf_{h\downarrow0} \frac{  \int_{
a}^{ b}\!\!\int_{ a}^{ b}|\rho_{h}(x-y)|^{2k_0+2 }\,dx\,dy}{ \int_{ a}^{
b}\!\!\int_{ a}^{ b}|\rho_{h}(x-y)|^{ 2k_0 }\,dx\,dy}>0 
 .\label{521ptq}
\ee
  Then
\be
\lim_{h\downarrow 0} {I_G(f,h;a,b)-(b-a)E f(\eta) \over
\sqrt{\mbox{Var
$I_G(f,h;a,b)$}}}
\stackrel{law}{=} N( 0,1).  \label{b3.1wwp}
\ee
\ec

\Proof 
We  write 
\be |\rho_{h}(x-y)|^{ 2k_0+ 1}=|\rho_{h}(x-y)|^{ k_0 }|\rho_{h}(x-y)|^{
k_0+1}\label{mm.56}
\ee
      and use the Schwarz Inequality to see that
\be
\frac{ \int_{ a}^{ b}\!\!\int_{ a}^{ b}|\rho_{h}(x-y)|^{2k_0+2}\,dx\,dy}{ \int_{
a}^{ b}\!\!\int_{ a}^{ b}|\rho_{h}(x-y)|^{2k_0+1}\,dx\,dy}\ge\frac{ \int_{
a}^{ b}\!\!\int_{ a}^{ b}|\rho_{h}(x-y)|^{2k_0+1 }\,dx\,dy}{ \int_{ a}^{
b}\!\!\int_{ a}^{ b}|\rho_{h}(x-y)|^{ 2k_0 }\,dx\,dy}.\label{521}
\ee
It follows from (\ref{521ptq}) that there exists a  $\de>0$ for which
\be
\liminf_{h\downarrow0} \frac{  \int_{
a}^{ b}\!\!\int_{ a}^{ b}|\rho_{h}(x-y)|^{2k_0+1 }\,dx\,dy}{ \int_{ a}^{
b}\!\!\int_{ a}^{ b}|\rho_{h}(x-y)|^{ 2k_0 }\,dx\,dy}=\de\label{521q}.
\ee
Consequently, for all $l>2k_0$
\be
\liminf_{h\downarrow0} \frac{  \int_{
a}^{ b}\!\!\int_{ a}^{ b}|\rho_{h}(x-y)|^{l }\,dx\,dy}{ \int_{ a}^{
b}\!\!\int_{ a}^{ b}|\rho_{h}(x-y)|^{ 2k_0 }\,dx\,dy}\ge\de^{l-2k_0}
.\label{521qq}
\ee
This shows that all  the integrals$ \int_{
a}^{ b}\!\!\int_{ a}^{ b}|\rho_{h}(x-y)|^{j}\,dx\,dy$ with  $2k_0\le j$ have
the same order of magnitude as $h$ decreases to zero. Therefore,
\begin{equation}
\(\int_a^b\!\!\int_a^b|\rho_h(x-y)|^j\,dx\,dy\)^{1/j}
=o\(\int_a^b\!\!\int_a^b
|\rho_h(x-y)|^{j+1}\,dx\,dy\)^{1/(j+1)}\label{bnee}
\end{equation}
for all $2k_0\le j$. This and (\ref{bneeq}) are all that is used in the
proof of Theorem \ref{CCLT}.
\qed

\medskip Lemma \ref{theo-weak2m} is stated for the $2k$-th Wick
power. It could just as well have been stated for the $2k$-th Hermite
polynomial. As such it gives just one term in the Hermite polynomial expansion
of $f\in  L^{ 2}(R^{ 1},\,d\mu )$. However, in some cases, depending on
$\si^2(h)$, this suffices to give the CLT for all $f$, as we show in the next
lemma.  

\bl\label{theo-GCLT2q} Let  $f\in L^{ 2}(R^{ 1},\,d\mu)$ be symmetric and let
\be k_{ 0}=\inf_{m\ge 1}\{m|a_{2m}\ne 0 \} \label{bmm.3.1q}
\ee
for $a_{2m}$ as given in (\ref{am.m1}). Suppose that 
\be
\lim_{h\downarrow0} \frac{  \int_{
a}^{ b}\!\!\int_{ a}^{ b}|\rho_{h}(x-y)|^{2k_0+2 }\,dx\,dy}{ \int_{ a}^{
b}\!\!\int_{ a}^{ b}|\rho_{h}(x-y)|^{ 2k_0 }\,dx\,dy}=0\label{521t}.
\ee
  Then
\begin{equation}
\mbox{Var
$I_G(f,h;a,b)$}\sim  a_{2k_{0}}^{2}\int_{ a}^{
b}\!\!\int_{ a}^{ b}|\rho_{h}(x-y)|^{ 2k_0 }\,dx\,dy\label{j.1v}
\end{equation}
and
\be
\lim_{h\downarrow 0} {I_G(f,h;a,b)-(b-a)E f(\eta) \over
\sqrt{\mbox{Var
$I_G(f,h;a,b)$}}} 
\stackrel{law}{=} N( 0,1).  \label{j.1q}
\ee
\el

\Proof It follows from (\ref{a18.1}),   and
(\ref{bmm.3.1q}), that
\bea && \int_{ a}^{ b} f \(\frac{G( x+h)- G( x) }{\si ( h)}\) \,dx
-(b-a)Ef(\eta)\label{m5.3}
\\ &&\qquad= a_{ 2k_0} \int_{ a}^{ b}
H_{2k_0}\(\frac{G( x+h)- G( x) }{\si ( h)}\) \,dx \nn\\&&\qquad \qquad
+\int_{ a}^{
b}\(\sum_{ m= k_0+1 }^{
\ff}a_{ 2 m}
\ H_{ 2m}\(  \frac{G( x+h)- G( x) }{\si ( h)} \) \)\,dx \nn\\
&&\qquad : =a_{2k_0} \int_{ a}^{ b}
H_{2k_0}\(\frac{G( x+h)- G( x) }{\si ( h)}\) \,dx + \int_{ a}^{ b}
W_h(x)\,dx.\nn
\eea
By   (\ref{ag2.3t}) and (\ref{1.6gq})
\bea
\mbox{Var } \(\int_{ a}^{ b}  W_h(x)\,dx\)&=&\sum_{ m=k_0+1}^{ \ff}a_{ 2m}^{
2}
\int_{ a}^{ b}\!\!\int_{ a}^{ b} |\rho_{h}(x-y)|^{ 2m}
\,dx\,dy\label{m5.3j}\\
&\le&\int_{ a}^{ b}\!\!\int_{ a}^{
b} |\rho_{h}(x-y)|^{2(k_0+1)}\,dx\,dy\(\sum_{ m=k_0+1}^{ \ff}a_{  2m}^{
2}\).\nn
\eea 
By (\ref{521t})
\begin{equation}
\mbox{Var } \(\int_{ a}^{ b}  W_h(x)\,dx\)=o\(\mbox{Var }
 I_G(f,h;a,b)\)\label{j.1w}
\end{equation}
and (\ref{j.1v}) follows.

By (\ref{m5.3})
\bea &&\lim_{h\to 0} \frac{\int_{ a}^{ b} f \(\frac{G( x+h)- G( x)
}{\si ( h)}\)
\,dx -(b-a)Ef(\eta)}{\sqrt{\mbox{Var
$I_G(f,h;a,b)$}}} \label{.m8}\\ &&\qquad=  \lim_{h\to
0}\frac{ \int_{ a}^{ b} H_{2k_0} \(\frac{G( x+h)- G( x) }{\si ( h)}\)
\,dx}{\sqrt{\mbox{Var
$I_G(f,h;a,b)$}}}\nn\\
& &\qquad \qquad +\lim_{h\to 0} \frac{\int_{ a}^{ b}
W_h(x)\,dx}{\sqrt{\mbox{Var
$I_G(f,h;a,b)$}}}\nn.
\eea
Using    (\ref{j.1w})   we see that
\be
\lim_{h\to 0}\,  \mbox{Var } \(\frac{\int_{ a}^{ b}
W_h(x)\,dx}{\sqrt{\mbox{Var
$I_G(f,h;a,b)$}}}\)= 0.
      \ee
Therefore, (\ref{j.1q}) follows from (\ref{.m8}), (\ref{j.1v}) and (\ref{m3.1}).
  \qed

\br\label{rem-3}{\rm   It  is easy to see that when (\ref{521t}) holds
\be {\mbox{Var
$I_G(f,h;a,b)$}}\sim\(E\(f(\eta)H_{2k_0}(\eta)\)\)^2\int_{ a}^{
b}\!\!\int_{ a}^{ b}|\rho_{h}(x-y)|^{ 2k_0 }\,dx\,dy 
\ee
and when (\ref{521ptq}) holds
\be   {\mbox{Var
$I_G(f,h;a,b)$}}\sim \sum_{m=k_0}^\ff \(E\(f(\eta)H_{2m}(\eta)\)\)^2\int_{ a}^{
b}\!\!\int_{ a}^{ b}|\rho_{h}(x-y)|^{ 2m }\,dx\,dy .
\ee
 }
\er

\section{Concave \boldmath{$\si^2$}}\label{sec-con}

        Using the fact that  $\rho_h$ is symmetric and setting $c=b-a$
we see that for all $k\in N$
\bea
\int_a^b\!\!\int_a^b
|\rho_h(x-y)|^k\,dx\,dy&=&\int_0^{c}\!\!\int_0^{c}
            |\rho_h(x-y)|^k\,dx\,dy\label{mm4.1}\\
&=&2\int_0^c|\rho_h(s)|^k(c-s)\,ds\nn .
\eea

       The function  $\si^2(h)$,   defined in (\ref{m1.2}), has the
properties that $\si^2(0)=0$, and  $\si^2(h)\not\equiv 0$.
Therefore,
        if  it is concave, it is also both increasing and strictly increasing
on
$[0,c_0]$, for some $c_0>0$. In what follows we assume that
$c=b-a\ge c_0$.

\bl\label{lem-7.1} When $\si^2(h)$ is concave on $[0,c]$,   for all
$0<h<<c$,  and $k\in N$
\bea
       {(c-h)\over 2^{k }} \int  ^{h}_0 | \si^2( h)- \si^2(
s)|^k\,ds&\le&\int_{ a}^{ b}\int_{ a}^{ b}|\phi_{h}(x-y)|^{
k}\,dx\,dy\label{7.53ss}\\ &  \le& 6c\,\(1+{1\over 2^k}\)\int
^{h}_0 | \si^2( h)- \si^2( s)|^k\,ds\nn.
\eea
\el

      The proof of Lemma
\ref{lem-7.1} uses the next lemma which is also used to give many
other properties of the integrals in (\ref{mm4.1}).

\bl \label{lem-7.2} When $\si^2(h)$ is concave, for all
$0<h<<c$,   and $k\in N$
\bea
\frac{1}{2^{k+1}}\int  ^{h}_0 | \si^2( h)- \si^2(
s)|^k\,ds&\le&\int_0^{ c}|\phi_{h}(s)|^{
k}\,ds\label{.m83}\\&\le&  3\(1+{1\over 2^k}\)\int ^{h}_0 |
\si^2( h)-
\si^2( s)|^k\,ds\nn.
\eea
\el

\Proof  It is useful to work with  $-\phi_{ h}(s)$ rather than
$\phi_{ h}(s)$. To avoid confusion we set  $\varphi_{ h}(s)=-\phi_{
h}(s)$. Obviously
$|\varphi_{ h}(s)|=|\phi_{ h}(s)|$. Using the fact that
$\si^2(s)$ is concave, we note that for
$0< s\le h$,
\bea
       \varphi_h(s) &=&{1\over
2}(\si^2(h+s)+\si^2(h-s)-2\si^2(s))\label{.m98w}\\ &\le&
\(\si^2(h)-\si^2(s)\).\nn
\eea Since $\si^2(s)$ is increasing, by  writing
\be
       \varphi_h(s)={1\over
2}(\si^2(h-s)-\si^2(s))+(\si^2(h+s)-\si^2(s))
\ee
         we see that
\begin{equation}
        \varphi_h(s)\ge 0\qquad \mbox{for $s\in[0,h/2]$}
.\label{.m98w7}
\end{equation} Let
\begin{equation} A_h:=\{ 0< s\le h\,|\,
\varphi_h(s)<0\}.\label{.m98w8}
\end{equation} Clearly $A_h\subset  (h/2,h]$. Furthermore, on
$A_h$, since $\si^2(s)$ is increasing
\bea |\varphi_h(s)|&=&{1\over
2}((\si^2(s)-\si^2(h-s))-(\si^2(h+s)-\si^2(s))\label{.m98w10}\\
&\leq&\frac{1}{2} (\si^2(s)-\si^2(h-s)).\nn
\eea Let
\begin{equation} B_h:=\{ 0< s\le h\,|\,0\leq
\varphi_h(s)\}\label{.m98w9}.
\end{equation} Then, by (\ref{.m98w}) and (\ref{.m98w10})
\bea  &&\int_0^h  |\varphi_h(s)|^k\,ds=\int_{B_h }
|\varphi_h(s)|^k\,ds+\int_{ A_h} |\varphi_h(s)|^k\,ds\\ &
&\qquad \le
\int_0^{h} |\si^2(h)-\si^2(s) |^k\,ds  +{1\over 2^k}\int_ {h/2}^h |
\si^2(s)-
\si^2(h-s)|^k\,ds\nn\\ &  &\qquad=  \int_0^{h}
|\si^2(h)-\si^2(s)|^k\,ds  +{1\over 2^k}\int  ^{h/2}_0 |
\si^2(s+h/2)-
\si^2(h/2-s)|^k\,ds\nn.
\eea

       Using the fact that $\si^2(s)$ is monotonically increasing, when
$0\leq s\leq h/2$, we have   $0\leq
\si^2(s+h/2)-
\si^2(h/2-s)\leq \si^2(h)- \si^2(h/2-s)$. Consequently,
\bea  &&\int  ^{h/2}_0 | \si^2(s+h/2)- \si^2(h/2-s)|^k\,ds\\  &
&\qquad \le\int  ^{h/2}_0 | \si^2( h)- \si^2(h/2-s)|^k\,ds\nn \\ &
&\qquad =\int  ^{h/2}_0 | \si^2( h)- \si^2( s)|^k\,ds,
\eea where the last step employs a simple change of variables. This
shows us that
\be
\int_0^h  |\phi_h(s)|^k\,ds=\int_0^h  |\varphi_h(s)|^k\,ds\le\(
1+{1\over 2^k}\)
\int  ^{h }_0 | \si^2( h)-
\si^2( s)|^k\,ds.\label{.m82}
\ee

Let $g$ be a convex increasing function with $g(0)=0$. Then, if
$a\ge b\ge0$,
$g(a-b)\le g(a)-g(b)$. Therefore, since $\si^2$ is concave and
increasing
\bea &&\int_ h^c g\(2|\phi_h(s)|\)\,ds\\ & &\qquad =\int_ h^c
g\(\(\si^2(s)-\si^2(s-h)\)-
\(\si^2(s+h)-\si^2(s)\)\)\,ds\nn\\& &\qquad \le \int_ h^c
g\(\si^2(s)-\si^2(s-h)\)\,ds-\int_ h^c g\(
\si^2(s+h)-\si^2(s)\)\,ds\nn\\ && \qquad = \int_ h^c
g\(\si^2(s)-\si^2(s-h)\)\,ds-\int_ {2h}^{c+h}
g\(\si^2(s)-\si^2(s-h)\)\,ds\nn\\ & &\qquad \le\int_ h^{2h}
g\(\si^2(s)-\si^2(s-h)\)\,ds\nn\\ & &\qquad =\int_ 0^h
g\(\si^2(s+h)-\si^2(s)\)\,ds\nn\\ & &\qquad \le 2\int_ 0^{h/2}
g\(\si^2(s+h)-\si^2(s)\)\,ds.\nn
\eea On the other hand, using (\ref{.m98w7})
\bea
\lefteqn{\int_ 0^{h/2}  g(2|\phi_h(s)|)\,ds\label{mm7.2}}\\
&=&\int_ 0^{h/2}
g\(\(\si^2(s+h)-\si^2(s)\)+\(\si^2(h-s)-\si^2(s)\)\)\,ds\nn\\
&\ge&\int_ 0^{h/2}
        g\(\si^2(s+h)-\si^2(s)\) \,ds.\nn
\eea Consequently,
\be
\int_ h^c g\(2|\phi_h(s)|\)\,ds\le 2\int_ 0^{h }  g(2|\phi_h(s)|)\,ds
       \ee and therefore
\be
\int_ 0^c g\(2|\phi_h(s)|\)\,ds\le 3\int_ 0^{h }
g(2|\phi_h(s)|)\,ds.\label{.m81}
       \ee Using  (\ref{.m81}) and (\ref{.m82}) with
$g(\cdot)=|\cdot|^k$ we get the upper bound in (\ref{.m83}).

To get the  lower bound in (\ref{.m83}) we note that
\bea
\int  ^{h}_0 |2 \phi_h(s))|^k\,ds\label{mm7.2a} &\ge &\int
^{h/2}_0 | 2\phi_h(s))|^k\,ds\\ &=& \int ^{h/2}_0
\((\si^2(h-s)-\si^2(s))+(\si^2(h+s)-\si^2(s))\)^k\,ds\nn\\&
\ge& \int ^{h/2}_0 |\si^2(h+s)-\si^2(s)|^k\,ds\nn\\ &\ge &\int
^{h/2}_0 |\si^2( h)- \si^2( s)|^k\,ds\nn
\eea
       which, since $\si^2( s)$ is increasing,  implies that
\be 2\int  ^{h}_0 |2 \phi_h(s))|^k\,ds\ge \int  ^{h }_0 | \si^2( h)-
\si^2( s)|^k\,ds.\label{.m61}
\ee
\qed

\medskip{\noindent}{\bf Proof of Lemma \ref{lem-7.1}} The upper
bound in (\ref{7.53ss}) follows immediately from Lemma
\ref{lem-7.2} and  (\ref{mm4.1}). Also, by (\ref{mm4.1}) and
(\ref{.m61})
\bea
\int_a^b\!\!\int_a^b
|2\phi_h(x-y)|^k\,dx\,dy&=&2\int_0^c|2\phi_h(s)|^k(c-s)\,ds\label{.m89}\\
&\ge&2(c-h)\int_0^h|2\phi_h(s)|^k\,ds\nn\\ &\ge& (c-h) \int  ^{h
}_0 |
\si^2( h)-
\si^2( s)|^k\,ds\nn.
\eea This gives the  lower bound in (\ref{7.53ss}). \qed

\medskip It is useful to record the following inequalities:

\bl\label{lem-6.3} When $\si^2(s)$ is concave on   $[0,c]$ it follows
that for some
$0<h<<c$,  and $k\in N$
\bea &&
\frac{1}{2c} \int_{ a}^{ b}\int_{ a}^{
b}|\rho_{h}(x-y)|^k\,dx\,dy\label{m1.5k} \\&&\qquad
\le\sup_{a\leq x\leq b }
\int_{ a}^{ b}|\rho_{h}(x-y)|^k\,dy\le  \frac{3}{c-h}\int_{ a}^{
b}\int_{ a}^{ b}|\rho_{h}(x-y)|^k\,dx\,dy \nn.
\eea
      In particular (\ref{a3.2}) holds if and only if
\begin{equation}
\liminf_{h\downarrow 0 }{ \int_{ a}^{ b}\int_{ a}^{
b}|\rho_{h}(x-y)|^{2 k}\,dx\,dy\over
       \int_{ a}^{ b}\int_{ a}^{ b}|\rho_{h}(x-y)|
      \,dx\,dy}>0.\label{3.2ee}
\end{equation}
\el

\Proof   For all $k\in N$
\bea
        \sup_{a\leq x\leq b }
\int_{ a}^{ b}|\rho_{h}(x-y)|^k\,dy&=& \sup_{a\leq x\leq b }
\int_{ a}^{ b}|\rho_{h}(y-x)|^k\,dy\label{mm4.56}\\ &=&
\sup_{a\leq x\leq b }
\int_{ a-x}^{ b-x}|\rho_{h}(s)|^k\,ds\nn.
\eea Using this and the fact that $\rho_h(s)=\rho_h(-s)$ we see that
\be
       \int_0^c|\rho_{h}(s)|^k\,ds\le   \sup_{a\leq x\leq b }
\int_{ a}^{ b}|\rho_{h}(x-y)|^k\,dy\le 2
\int_0^c|\rho_{h}(s)|^k\,ds.\label{mm4.1.1}
\ee

Using (\ref{mm4.1}) and (\ref{mm4.1.1}) we see that
       \bea
        \sup_{a\leq x\leq b }
\int_{ a}^{
b}|\rho_{h}(x-y)|^k\,dy&\ge&\frac1c\int_0^c|\rho_h(s)|^k(c-s)\,ds\label{mm.5.1}\\
      &=&\frac{1}{2c}\int_a^b\!\!\int_a^b
|\rho_h(x-y)|^k\,dx\,dy.\nn
\eea This gives the first inequality is given in (\ref{m1.5k}).
       For the second inequality we see that
        by   (\ref{mm4.1.1}), (\ref{.m81}) and (\ref{mm4.1})
       \bea
          \sup_{a\leq x\leq b }
\int_{ a}^{ b}|\rho_{h}(x-y)|^k\,dy\label{m1.5kff}&\le& 2\int_{
0}^{ c}|\rho_{h}(s)|^k\,ds\\ &\le& \frac{6}{c-h}\int_{ 0}^{
h}|\rho_{h}(s)|^k(c-s)\,ds\nn\\ &\le& \frac{6}{c-h}\int_{ 0}^{
c}|\rho_{h}(s)|^k(c-s)\,ds\nn\\ &=&  \frac{3}{c-h} \int_{ a}^{
b}\int_{ a}^{ b}|\rho_{h}(x-y)|^k\,dx\,dy \nn .
\eea The rest of the lemma is obvious.
\qed

\bl  \label{lem-7.2f}  When $\si^2(h)$ is concave on $[0,c]$,
(\ref{bneed.1aa})   and (\ref{am1.6}) hold   for all $k\in N$. In addition
\begin{equation}
\(\int_a^b\!\!\int_a^b|\rho_h(x-y)|^k\,dx\,dy\)^{1/k}
=o\(\int_a^b\!\!\int_a^b
|\rho_h(x-y)|^{k+1}\,dx\,dy\)^{1/(k+1)} \label{bneed.2rq}
\end{equation}
for all $k\in N$, so
(\ref{bneed.2r})  also holds.
\el

\Proof
       Using (\ref{mm4.1.1})   and (\ref{.m81}) we see that
\bea
       \sup_{a\leq x\leq b }
\int_{ a}^{ b}|\rho_{h}(x-y)|^k\,dy&\le& 2\int_0^c|\rho_h(s)|^k
\,ds\label{.m69v}\\ &\le& 6\int_0^h|\rho_h(s)|^k \,ds,\nn
\eea 
which gives (\ref{am1.6}). (Here we use the simple observation
that if
\be
\int_ h^c   |\phi_h(s)|^k\,ds\le 2\int_ 0^{h }   |\phi_h(s)|^k\,ds
       \ee then
\be
\int_ h^c   |\rho_h(s)|^k\,ds\le 2\int_ 0^{h }   |\rho_h(s)|^k\,ds.
       \ee We  continue to pass between relations for $\phi$ and
$\rho$ in this way without further comment.)

The condition in  (\ref{bneed.1aa}) follows from (\ref{m1.5k}).

To obtain (\ref{bneed.2rq}) we note that for
$k<m\in N$, by (\ref{7.53ss}) used twice
\bea &&\(\int_{ a}^{ b}\int_{ a}^{ b}|\phi_{h}(x,y)|^{
k}\,dx\,dy\)^{1/k}\label{mm.56jr}\\ &&\qquad \le C_1\(\int
^{h}_0
  |  \si^2( h)-
\si^2( s)  |^k\,ds\)^{1/k}\nn\\  &&\qquad \le
C_1h^{1/k-1/m}\(\int ^{h}_0  \ |  \si^2( h)-
\si^2( s |^m\,ds\)^{1/m}\nn \\  &&\qquad \le
C_1h^{1/k-1/m}\(C_2\int_{ a}^{ b}\int_{ a}^{ b}|\phi_{h}(x,y)|^{
m}\,dx\,dy\)^{1/m}\nn,
\eea where $C_1$ and $C_2$ are finite constants that only depend
on
$c=b-a$ for all $h<<c$.\qed

\begin{lemma}\label{lem-regpos}Let  $\si^2(h)$ be concave and
regularly varying with  index $\ga>0$. Then for all $k\in N$
\begin{equation}
\int_{ a}^{ b}\int_{ a}^{ b}|\rho_{h}(x-y)|^{
k}\,dx\,dy\approx h.\label{j.7}
\end{equation}
\end{lemma}

\Proof
It is clear from (\ref{7.53ss}) that
\be
       \int_{ a}^{ b}\int_{ a}^{ b}|\rho_{h}(x-y)|^{ k}\,dx\,dy   \le\,\,
6c\(1+{1\over 2^k}\) h\label{523}.
\ee
Also by
Lemma
\ref{lem-7.1}
\bea
\int_{ a}^{ b}\int_{ a}^{ b}|\rho_{h}(x-y)| \,dx\,dy&\ge&
\frac{c-h}{2} \int ^{h}_0
\bigg|1-\frac{\si^2(s)}{\si^2(h)}\bigg| \,ds\label{.m7}\\ &\ge&  \frac{c-h}{2}
\int ^{h/2}_0
\bigg|1-\frac{\si^2(h/2)}{\si^2(h)}\bigg| \,ds\nn.
\eea
 When $\si^2(h)$ is regularly varying at zero with index
$\ga>0$
\be
\lim_{h\to 0}\bigg|1-\frac{\si^2(h/2)}{\si^2(h)}\bigg|=1-{1 \over 2^{\ga}}. 
\ee
Using this in (\ref{.m7}) we get the lower bound in (\ref{j.7}).
\qed

\medskip
In preparation  for the next lemma we point out that when $\si^2(s)$ is
concave and regularly varying with index  $\ga\ge 0$
\be
\lim_{s\to 0}\frac{s\,\frac{d}{ds}\( \si^2(s)
\)}{\si^2(s)}=\ga.\label{611}
\ee This follows from the Monotone Density Theorem
\cite[Theorem1.7.2b]{BGT}, (see also
\cite[page 596]{book}),   since the  derivative of $\si^2(s)$ is
decreasing.
\bl\label{lem-4.5} Let $\si^2(s)$ be concave on $[0,c]$.

\begin{itemize}
\item[(1)] If $s\frac{d}{ds}\si^2(s)$ is increasing on $[0,h]$, then
for some $0<h< <c$   and all $k\ge 1$,
\be
       \int_{ a}^{ b}\int_{ a}^{ b}|\rho_{h}(x-y)|^{
k}\,dx\,dy\label{7.53dad}    \le C_{c,k}
\(\frac{h
       \frac{d}{dh}\si (h)}{\si(h)}  \)^k h,
\ee where  $ C_{c,k}<\ff$ depends only on $c$ and $k$.
\item[(2)] If $\si^2(s)$ is slowly varying at zero  then   for all $k\ge 1$,
\end{itemize}
\be
      \limsup_{h\downarrow 0} {\int_{a}^{ b}\int_{ a}^{
b}|\rho_{h}(x-y)|^{ k}\,dx\,dy\over h}\label{7.53pop}=0,
\ee and
\begin{equation}
\limsup _{h\downarrow 0 }{ \int_{ a}^{ b}\int_{ a}^{
b}|\rho_{h}(x-y)|^{ k}\,dx\,dy\over
       \int_{ a}^{ b}
\int_{ a}^{ b}|\rho_{h}(x-y)|^{k-1} \,dx\,dy}=0.\label{602}
\end{equation}
\el

\Proof By Lemma \ref{lem-7.1}
\be
       \int_{ a}^{ b}\int_{ a}^{ b}|\rho_{h}(x-y)|^{
k}\,dx\,dy\label{7.53sasq}  \le 6c\(1+{1\over 2^k}\) \int  ^{h}_0
\bigg| 1- \frac{\si^2( s)}{\si^2( h)}
\bigg|^k\,ds.
\ee
       Using integration by parts we see that
\bea
\int  ^{h}_0 \bigg| 1- \frac{\si^2( s)}{\si^2( h)}
\bigg|^k\,ds&=&\frac{k}{\si^2( h)}\int  ^{h}_0
\bigg| 1- \frac{\si^2( s)}{\si^2( h)} \bigg|^{k-1} s\,\frac{d}{ds}\(
\si^2(s)
\) \,ds  \,\label{.m63q}\\ & \le &\frac{kh}{\si^2( h)}\frac{d}{dh}\(
\si^2(h) \)\int  ^{h}_0
\bigg| 1- \frac{\si^2( s)}{\si^2( h)} \bigg|^{k-1}
       \,  \,ds\nn,
\eea where at the last step we use the fact that
$s\frac{d}{ds}\si^2(s)$ is increasing on $[0,h]$.  Since the last
integral in (\ref{.m63q}) is equal to
$h$ when $k=1$ we get (\ref{7.53dad}).

To obtain (\ref{7.53pop}) we use  the first line of (\ref{.m63q}) to
get
\bea
\int  ^{h}_0 \bigg| 1- \frac{\si^2( s)}{\si^2( h)}
\bigg|^k\,ds&= &\frac{k}{\si^2( h)}\int  ^{h}_0
\bigg| 1- \frac{\si^2( s)}{\si^2( h)} \bigg|^{k-1} s\,\frac{d}{ds}\(
\si^2(s)
\) \,ds  \,\label{.m63aq}\\ & \le &k\int  ^{h}_0
\bigg| 1- \frac{\si^2( s)}{\si^2( h)}
\bigg|^{k-1}\frac{s\,\frac{d}{ds}\( \si^2(s)
\)}{\si^2(s)}
       \,  \,ds\nn.
\eea Consequently, it follows from (\ref{611}),   with $\ga=0$, that
\be
\limsup_{h\downarrow 0}\frac{\int  ^{h}_0 \bigg| 1- \frac{\si^2(
s)}{\si^2( h)}
\bigg|^k\,ds}{\int  ^{h}_0 \bigg| 1- \frac{\si^2( s)}{\si^2( h)}
\bigg|^{k-1}\,ds}=0 .\label{601}
\ee Iterating this and using (\ref{7.53sasq}) we get (\ref{7.53pop}).

The statement in (\ref{602}) follows from (\ref{601}) and Lemma
\ref{lem-7.1}. \qed

     \medskip By the first line of (\ref{.m63q})
\be
\int  ^{h}_0 \bigg| 1- \frac{\si^2( s)}{\si^2( h)} \bigg|
\,ds=\frac{1}{\si^2( h)}\int  ^{h}_0 s\,\frac{d}{ds}\( \si^2(s) \)
\,ds.
\ee
     By (\ref{611}) when $\si^2(s)$ is concave and regularly varying
with index $ \ga\ge 0$,
$s\,\frac{d}{ds}\(
\si^2(s)
\)
$ is regularly varying with   index $ \ga\ge 0$.  Using this we see
that
\bea
\frac{1}{\si^2( h)}\int  ^{h}_0 s\,\frac{d}{ds}\( \si^2(s) \)
\,ds&\sim&
\frac{1}{1+\ga}
\frac{h^2\frac{d}{dh}\(
\si^2(h)
\)}{\si^2(h)}\\ &  = &\frac{2}{1+\ga} \frac{h^2\frac{d}{dh}\( \si
(h)
\)}{\si (h)}.\quad
\eea Therefore, it follows from (\ref{7.53ss})
\be {  \int_{ a}^{ b}\int_{ a}^{ b}|\rho_{h}(x-y)| \,dx\,dy \over
h}\approx \frac{h
\frac{d}{dh}\(
\si (h)
\)}{\si (h)}\label{mm.58}.
       \ee When $\si^2(h)$ is concave, the right-hand side of
(\ref{mm.58}) goes to
$\ga$ as  $h$ decreases to zero. For $\ga>0$, this restates a
property given in Lemma
\ref{lem-regpos}. However, when $\ga=0$ this is a refinement of
(\ref{7.53pop}).

\bl \label{lem-variance}  When $\si^2(h)$ is concave   on $[0,c]$
and
     regularly varying with index $\ga\ge 0$ and
$s\frac{d}{ds}\si^2(s)$ is increasing  for some
$0<h< <c$  and all $k\ge N$.
\be
       \int_{ a}^{ b}\!\!\int_{ a}^{ b}|\rho_{h}(x-y)|^{ k}\,dx\,dy
\approx\(
\frac{h \frac{d}{dh}
\si (h)
  }{\si(h)} \)^{k}h \label{mm.57q}.
\ee
\el

\Proof By (\ref{mm.56jr}) 
\be
C_1 h^{1-1/k}  \(\int_{ a}^{ b}\!\!\int_{ a}^{ b}|\rho_{h}(x-y)|^{
k}\,dx\,dy\)^{1/k}
\ge
  \int_{ a}^{ b}\!\!\int_{ a}^{ b}|\rho_{h}(x-y)|
\,dx\,dy \label{511}
\ee or, equivalently
\be
       \int_{ a}^{ b}\!\!\int_{ a}^{ b}|\rho_{h}(x-y)|^{ k}\,dx\,dy \ge
C_1^{-k}\(\frac{ \int_{ a}^{ b}\!\!\int_{ a}^{ b}|\rho_{h}(x-y)|
\,dx\,dy}{h}\)^{k }h.\label{mm.57}
\ee
       where $C_1>0$   depends only on
$c=b-a$ for all  $h<<c$. Using this, (\ref{mm.58}) and (\ref{7.53dad})
completes the proof. \qed

\medskip It follows from Lemma \ref{lem-6.3}
        that when
\begin{equation}
\lim _{h\downarrow 0 }{ \int_{ a}^{ b}\int_{ a}^{
b}|\rho_{h}(x,y)|^{ k}\,dx\,dy\over
\int_{ a}^{ b}\int_{ a}^{ b}|\rho_{h}(x,y)| \,dx\,dy}=0\label{3.2wy}
\end{equation} the condition in (\ref{a3.2}) fails. Therefore, by  the
first line of (\ref{.m63q}) and Lemma \ref{lem-7.1}, if
$\si^2(h)$ is concave and
\be
\lim_{h\downarrow 0}h\frac{d}{dh}\(\log\si^2(h)
\)=\lim_{h\downarrow 0}\frac{h}{\si^2(h)} \frac{d}{dh}
\(  \si^2(h)
\)=0\label{.m74}
\ee (\ref{a3.2}) fails.

\bl\label{lem-6.4}Assume that  $\si^2(s)$ is concave on $[0,h]$ for
some
$0<h<<c$. Write
\be
       \si^2(s)=\exp\(f(\log 1/s)\)
\ee If $\lim_{x\to \ff}f'(x)=0$, (\ref{a3.2}) fails.

\el

\Proof This is simple since
\be s \frac{d}{ds} \(\log \si^2(s) \)=s \frac{d}{ds}\(f(\log 1/s)\)=-f'
(\log 1/s).\label{.m75}
\ee The assertion follows from (\ref{.m74}). \qed

\section{\boldmath{ $\si^2(h)=h^r$, $0<r\le 2$} }\label{sec-sigsq}

In these cases we can  find precise asymptotic limits at zero of the
double integral in (\ref{7.53ss}) and thus obtain a precise value for
\mbox{Var
$I_G(f,h;a,b)$}. 
       We begin with the following estimates:

\bl\label{lem-ex1} Let  $\si^2(h)=h^r$, $0<r\le 2$.
\newline
        When $(2-r)k<1$
\be
           \int_0^{c} |\phi_h(s)|^{ k}\,ds \sim { r^{ k}| r-1|^{ k}c^{
(r-2)k+1}\over 2^k (r-2)k+1}h^{2k}\qquad \mbox{at zero}.
\label{61.1}
\ee When $r=1$
\be
           \int_0^{c} |\phi_h(s)|^k \,ds ={h^{k+1}\over k+1}
\label{61.3a}.
\ee When $(2-r)k=1$, $k\ge 2$
\be
           \int_0^{c} |\phi_h(s)|^{ k}\,ds \sim \Big|{ r( r-1)\over
2}\Big|^{k}h^{2k }\log 1/h\qquad \mbox{at zero}.
\label{61.3}
\ee
       If $(2-r)k>1$
\be
          \int_0^{c} |\phi_h(s)|^{ k}\,ds\sim h^{rk+1 }\int_0^{\ff}
|\phi_1(s)|^{ k}\,ds
\qquad \mbox{at zero}.
\label{61.2}
\ee
\el

\Proof  The equality in (\ref{61.3a}) is a trivial direct computation.
We proceed to the others. By a  simple change of variables we have
\begin{eqnarray} &&
    \int_0^{c} |2\phi_h(s)|^{ k}\,ds\label{61.4}\\ &&\qquad
=\int_0^{c}
\Big| |s+h|^{ r} +|s-h|^{ r}-2|s|^{ r}   \Big |^{ k}\,ds\nonumber\\
&&\qquad =h^{rk+1 }\int_0^{c/h} \Big| |s+1|^{ r} +|s-1|^{ r}-2|s|^{
r}   \Big |^{ k}\,ds\nonumber\\ &&\qquad =h^{rk+1 }\int_0^{c/h}
|2\phi_1(s)|^{ k}\,ds.\nonumber
\end{eqnarray} We write
\begin{eqnarray} &&
\int_0^{c/h} \Big| |s+1|^{ r} +|s-1|^{ r}-2|s|^{ r}   \Big |^{
k}\,ds\label{61.5}\\ &&\qquad =\int_0^{2} \Big| |s+1|^{ r} +|s-1|^{
r}-2|s|^{ r}   \Big |^{ k}\,ds\nonumber\\&&\qquad  \quad
+\int_2^{c/h}
\Big| |s+1|^{ r} +|s-1|^{ r}-2|s|^{ r}   \Big |^{ k}\,ds.\nonumber
\end{eqnarray} The first integral is   a finite number. For the second
integral we have, for
$h<c/2$
\begin{eqnarray} &&
\int_2^{c/h} \Big| |s+1|^{ r} +|s-1|^{ r}-2|s|^{ r}   \Big |^{
k}\,ds\label{61.6}\\ &&\qquad =\int_2^{c/h} s^{ rk}\Big| |1+s^{
-1}|^{ r} +|1-s^{ -1}|^{ r}-2   \Big |^{ k}\,ds\nonumber\\ &&\qquad
=\int_2^{c/h} s^{ rk}\Big| r( r-1)s^{ -2}+O( s^{ -3})\Big|^{
k}\,ds\nonumber\\ &&\qquad =\int_2^{c/h} s^{ (r-2)k}\Big| r(
r-1)+O( s^{ -1})\Big|^{ k}\,ds.\nonumber
\end{eqnarray} Using (\ref{61.4})--(\ref{61.6}) we get (\ref{61.1}),
(\ref{61.3}) and (\ref{61.2}). For (\ref{61.1}) and (\ref{61.3}), to do
the integration, it is helpful to note that when
$ (2-r)k<1$,  $ (r-2)k>-1$ and, obviously, when  $ (2-r)k=1$,   $
(r-2)k=-1$. For (\ref{61.2}) we have $ (r-2)k<-1$ so that the last
integral  in (\ref{61.6}),  and hence in (\ref{61.5}), is finite.

\qed

\medskip We now consider the integral in (\ref{mm4.1}).

\bl\label{lem-ex2}  Let  $\si^2(h)=h^r$, $0<r\le 2$. \newline When
$(2-r)k<1$
\be
\int_a^b\!\!\int_a^b |\phi_h(x-y)|^k\,dx\,dy \sim {2 r^{ k}| r-1|^{
k}c^{ (r-2)k+2}\over 2^k((r-2)k+1)((r-2)k+2)}h^{2k}\qquad
\mbox{at zero}.
\label{61.11}
\ee When $r=1$
\be
          \int_a^b\!\!\int_a^b |\phi_h(x-y)|^{ k} \,dx\,dy \sim
2c\,{h^{k+1}\over k+1}\qquad \mbox{at zero}.
\label{61.11ab}
\ee When $(2-r)k=1$, $k\ge2$
\be
          \int_a^b\!\!\int_a^b |\phi_h(x-y)|^k\,dx\,dy \sim 2c\Big|{ r(
r-1)\over 2}\Big|^{k}h^{2k }\log 1/h\qquad \mbox{at zero}.
\label{61.11a}
\ee When $(2-r)k>1$
\be
         \int_a^b\!\!\int_a^b |\phi_h(x-y)|^k\,dx\,dy\sim 2ch^{rk+1
}\int_0^{\ff} |\phi_1(s)|^{ k}\,ds \qquad \mbox{at zero}.
\label{61.12}
\ee
\el

\Proof     By (\ref{mm4.1}) it suffices to consider
\begin{equation}
2\int_0^c|\phi_h(s)|^k(c-s)\,ds=2c\int_0^c|\phi_h(s)|^k\,ds-2
\int_0^c|\phi_h(s)|^k\,s\,ds.\label{61.10a}
\end{equation}
As in the proof of Lemma
\ref{lem-ex1},   (\ref{61.11ab}) is a trivial direct computation. We
consider  the others. The first integral on the right-hand side of
(\ref{61.10a})  is handled by Lemma
\ref{lem-ex1} and obviously gives the results in Lemma
\ref{lem-ex1} multiplied by $2c$. For the last integral, by a change
of variables,  we have
\begin{eqnarray} &&
\int_0^{c} |2\phi_h(s)|^{ k}\,s\,ds\label{61.14}\\ &&\qquad
=\int_0^{c} \Big| |s+h|^{ r} +|s-h|^{ r}-2|s|^{ r}   \Big |^{
k}\,s\,ds\nonumber\\ &&\qquad    =h^{rk+2 }\int_0^{c/h} \Big|
|s+1|^{ r} +|s-1|^{ r}-2|s|^{ r}   \Big |^{ k}\,s\,ds\nonumber\\
&&\qquad =h^{rk+2 }\int_0^{c/h} |2\phi_1(s)|^{
k}\,s\,ds.\nonumber
\end{eqnarray}
In the case of (\ref{61.12})  as in (\ref{61.6}) we can bound
(\ref{61.14}) by
\begin{equation}
Ch^{rk+2 }\int_{ 2}^{ c/h} { 1\over s^{ ( 2-r)k-1}}\,ds.\label{j.8}
\end{equation}
If $( 2-r)k>2$ the integral is bounded whereas if $( 2-r)k=2$ the
integral $\approx \log 1/h$.
Thus the last integral in (\ref{61.14}) contributes nothing to the
asymptotic estimate of (\ref{61.10a}) at zero in these cases. When
   $1<( 2-r)k<2$ we see that  (\ref{j.8}) is  equal to $Ch^{rk+2 }h^{ (
2-r)k-2}=Ch^{rk+( 2-r)k}=o(h^{rk+1})$ since $1<( 2-r)k$. Hence
the last integral in (\ref{61.14}) contributes nothing to the asymptotic
estimate of (\ref{61.10a}) at zero in this case as well.  

In the cases of (\ref{61.11}) and (\ref{61.11a})
       we compute the integral in   (\ref{61.14}) using   (\ref{61.5}) and
(\ref{61.6}). We see that it  contributes nothing to the asymptotic estimate 
at zero in (\ref{61.11a}) but it does enter into the estimates in
(\ref{61.11}). \qed

\medskip We write the estimates in Lemma \ref{lem-ex2} in
different forms that are useful to us.

\begin{corollary}\label{cor-7}  Let  $\si^2(h)=h^r$, $0<r\le 2$.
\[
\(\int_a^b\!\!\int_a^b |\phi_h(x-y)|^k\,dx\,dy \)^{ 1/k}
\sim\left\{\begin{array}{lr} D_{1,k}\,h^{ 2}&\mbox{ $(2-r)k<1$}\\
\({2c\over k+1}\)^{1/k}\,h^{1+1/k}&\mbox{$r=1$ }\\ D_{2,k}\,h^{
2}(\log 1/h)^{ 1/k}&\mbox{ $(2-r)k=1$, $k\ge 2$}\\ D_{3,k}\,h^{
r+1/k}&\mbox{
$(2-r)k>1$}
\end{array}\right.
\]
       Also
\[
       \int_a^b\!\!\int_a^b |\rho_h(x-y)|^k\,dx\,dy
\sim\left\{\begin{array}{lr} D_{4,k}\,h^{ (2-r)k}&\mbox{$(2-r)k<1$}\\
{2c\over k+1}\,h &\mbox{$r=1$}\\ D_{5,k}\,h (\log 1/h)
&\mbox{$(2-r)k=1$. $k\ge 2$}\\ D_{6,k}\,h &\mbox{$(2-r)k>1$}
\end{array}\right.
\] Here $D_{j,k}=D_{j,k}(r ,c)$, $j=1,\ldots, 6$, do not depend on $h$.
(They can be obtained from Lemma \ref{lem-ex2}.)
\end{corollary}

\section{Proofs of Theorems \ref{theo-GCLT} and
\ref{theo-GCLT2} and Tables 1 and 2}\label{sec-proofs}

\noindent{\bf Proof of Theorem \ref{theo-GCLT} } All we need to do is verify
that the hypotheses of Theorem \ref{BCLT} are satisfied. When $\si^2(h)$ is
concave we show this in Lemma \ref{lem-7.2f}. It remains to consider 
$\si^2(h)=h^r$, $1<r\le 3/2$. As we show in (\ref{mm4.1.1})
\be
   \sup_{a\leq x\leq b }
\int_{ a}^{ b}|\rho_{h}(x-y)|^k\,dy\le 2
\int_0^c|\rho_{h}(s)|^k\,ds.\label{mm4.1.1p}
\ee
and, as we show in (\ref{mm4.1})
\be 
\int_a^b\!\!\int_a^b
|\rho_h(x-y)|^k\,dx\,dy = 2\int_0^c|\rho_h(s)|^k(c-s)\,ds.\label{mm4.1.1pp}
\ee
One can see from Lemmas \ref{lem-ex1} and \ref{lem-ex2} that the
right-hand sides of (\ref{mm4.1.1p}) and (\ref{mm4.1.1pp}) have the same
asymptotic behavior at zero for all $\si^2(h)=h^r$,
$0<r<2$. Thus we have (\ref{bneed.1aa}) when $\si^2(h)=h^r$, $1<r<2$.

We now show that  when  $\si^2(h)=h^r$, $1<r\le 3/2$ 
\begin{equation}
\(\int_a^b\!\!\int_a^b|\rho_h(x-y)|^j\,dx\,dy\)^{1/j}
=o\(\int_a^b\!\!\int_a^b
|\rho_h(x-y)|^{2k_0}\,dx\,dy\)^{1/(j+1)}\label{ab.2r}
\end{equation}
for all $j\in N$ which, of course, implies (\ref{bneed.2r}). By
Corollary \ref{cor-7}
\be
\int_a^b\!\!\int_a^b
|\rho_h(x-y)| \,dx\,dy\sim D_{4,1}\,h^{(2-r)}
\ee and when $j>2$
\be
\int_a^b\!\!\int_a^b
|\rho_h(x-y)|^j \,dx\,dy\sim D_{6,j}\,h \label{.m69}
\ee
and
\be
       \int_a^b\!\!\int_a^b |\rho_h(x-y)|^2\,dx\,dy
\sim\left\{\begin{array}{lr} D_{5,2}h (\log 1/h)
&\mbox{$r=3/2$ }\\ D_{6,2}h &\mbox{$1<r<3/2$}
\end{array}\right.\label{m.m9}.
\ee  
When $1<r< 3/2$, $2-r>1/2$ and (\ref{.m69}) holds for all $j\ge 2$. Thus we
get  (\ref{ab.2r}). When $r= 3/2$, $2-r=1/2$ but we get the extra $\log 1/h$
term in (\ref{m.m9}) so  we
get  (\ref{ab.2r}) in this case as well.\qed

\medskip\noindent{\bf Proof of Theorem \ref{theo-GCLT2} }    We
show  that the hypotheses of Corollary
\ref{CCLT} are satisfied. We already showed, in the the proof of 
Theorem \ref{theo-GCLT}, that (\ref{bneed.1aa}) holds for $\si^2(h)=h^r$,
$1<r<2$ so, in particular it holds for $3/2<r\le 2-1/(2k_0)$. Suppose
$r= 2-1/(2k_0)$. Then by Corollary \ref{cor-7}
\be
\(\int_a^b\!\!\int_a^b
|\rho_h(x-y)|^j \,dx\,dy\)^{1/j}\sim (D_{4,1})^{1/j}\,h^{(2-r)}\label{finishing}
\ee
and 
\be
\(\int_a^b\!\!\int_a^b
|\rho_h(x-y)|^{2k_0} \,dx\,dy\)^{1/(2k_0)}\sim (D_{5,2}h (\log
1/h))^{1/(2k_0)}.
\ee
Since, in this case $2-r=1/(2k_0)$ we see that (\ref{bneeq}) holds. Also, by
Corollary \ref{cor-7}, when $j>2k_0$, $(2-r)j>1$, and 
\be
\int_a^b\!\!\int_a^b |\rho_h(x-y)|^j\,dx\,dy\sim D_{6,k}h.\label{.m.m}
\ee
Thus (\ref{521ptq}) is also satisifed.

When  $3/2<r<2-1/(2k_0)$ it follows from Corollary \ref{cor-7} that 
\be
\(\int_a^b\!\!\int_a^b
|\rho_h(x-y)|^{2k_0} \,dx\,dy\)^{1/(2k_0)}\sim  (D_{6,2}h )^{1/(2k_0)} 
\ee
and for $j<2k_0$ 
\[
       \(\int_a^b\!\!\int_a^b |\rho_h(x-y)|^j\,dx\,dy\)^{1/j}
\sim\left\{\begin{array}{lr} (D_{4,j})^{1/j}\,h^{ (2-r) }&\mbox{$(2-r)j<1$}\\ 
(D_{5,k} \, h (\log 1/h))^{1/j} &\mbox{$(2-r)j=1$ }\\( D_{6,k}\,h)^{1/j}
&\mbox{$(2-r)j>1$}
\end{array}\right.
\]
Since, in this case both  $2-r>1/(2k_0)$ and $1/j>1/(2k_0)$,  
(\ref{bneeq}) holds.  When $j\ge 2k_0$  we are in the same situation as in
(\ref{.m.m}) so (\ref{521ptq}) is also satisifed.\qed

\medskip\noindent{\bf Explanation of how the entries in Table 
1 are obtained: } Entries $(1)-(4)$ are given in Corollary \ref{cor-7}. Entry (5)
is given in Lemma \ref{lem-regpos}. Entries (6) and (7) follow from Lemma
\ref{lem-variance}. The constants in (\ref{constants}) are taken from
Lemma \ref{lem-ex2}. \qed

\medskip\noindent{\bf Explanation of how the entries in Table 
2 are obtained: } As we point out just before Corollary \ref{theo-1}, 
$k_0=1$, and 
$a_2=E(|\eta|^p|\eta^2-1|)/\sqrt2>0$. The variance $\Phi(h)$ is given in
(\ref{j.5}) and we get the asymptotic estimates for $
\int_{ a}^{ b}\!\!\int_{ a}^{ b} (\rho_{h}(x-y))^{2 k}
\,dx\,dy$ from Table 1 for (1)-(3) and from Lemma \ref{lem-variance} for (5).
Recall Remark
\ref{rem-3}. In (1) and (5), (\ref{521t}) holds so $\Phi(h)$ is the 
single term $a_{ 2}^{ 2} \int_{ a}^{ b}\!\!\int_{ a}^{ b}
(\rho_{h}(x-y))^{2 }
\,dx\,dy$. In (2)   we get the infinite series. Example  (3) is simply (2) with
the integral evaluated. For  (4) we see by Lemma \ref{lem-regpos} that
(\ref{521ptq}) holds. Since the variance contains an infinite number of terms
we also need to use (\ref{523}) to get the estimate for $\Phi(h)$.  \qed

\medskip When $k_0\ge 2$ and $\si^2(h)=h^r$, $3/2<r\le 2-1/(2k_0)$ we
can also get precise asymptotic estimates for the denominator in (\ref{j.1}).
We leave this to the interested reader.

\def\noopsort#1{} \def\printfirst#1#2{#1}
\def\singleletter#1{#1}
               \def\switchargs#1#2{#2#1}
\def\bibsameauth{\leavevmode\vrule height .1ex
               depth 0pt width 2.3em\relax\,}
\makeatletter
\renewcommand{\@biblabel}[1]{\hfill#1.}\makeatother
\newcommand{\bysame}{\leavevmode\hbox to3em{\hrulefill}\,}

\end{document}